\documentclass[12pt]{amsart}
\usepackage[utf8]{inputenc}
\usepackage[margin=1in]{geometry}

\usepackage{changes}
\usepackage{comment}
\usetikzlibrary{patterns}
\usepackage{accents}

\usepackage{fullpage}
\usepackage{amssymb,amsfonts}
\usepackage[all,arc]{xy}
\usepackage{enumerate}
\usepackage{enumitem}
\usepackage{mathrsfs}
\usepackage{mathtools}
\usepackage{tikz}
\usepackage{tikz-cd}
\usepackage{graphicx}
\usepackage[font={small}]{caption}
\usepackage{float}
\usepackage{subcaption}
\usepackage{hyperref}
\usepackage{rotating}
\usepackage{changes}
\usepackage{standalone}
\usepackage{changes}
\usepackage{relsize}

\usepackage{pgfplots}
\pgfplotsset{compat=1.11}
\pgfplotsset{ticks=none}

\DeclareMathOperator{\Hom}{Hom}

\DeclareMathOperator{\F}{\mathbb{F}}
\DeclareMathOperator{\Z}{\mathbb{Z}}
\DeclareMathOperator{\Q}{\mathbb{Q}}
\DeclareMathOperator{\R}{\mathbb{R}}

\DeclareMathOperator{\Hbb}{\mathbb{H}}

\DeclareMathOperator{\Spinc}{Spin^c}
\DeclareMathOperator{\spinc}{spin^c}
\DeclareMathOperator{\Char}{Char}

\newcommand{\til}[1]{\widetilde{#1}}

\newcommand{\Zhat}{\widehat{Z}}
\newcommand{\normalization}{\Delta}
\newcommand{\tnormalization}{\Theta}
\newcommand{\h}[1]{\widehat{#1}}
\newcommand{\V}{\mathcal{V}}
\newcommand{\sgn}{\operatorname{sgn}}
\newcommand{\lr}[1]{\vert {#1} \vert}

\newcommand{\lattice}{%
\tikz[baseline]{\draw (0,0) --(0,1ex);
\draw (.5ex,0) --(.5ex,1ex);
\draw (-.4ex,.7ex) --(.9ex,.7ex);
\draw (-.4ex,.3ex) --(.9ex,.3ex);
}%
}

\newcommand{\ZhatL}{\Zhat^{\lattice}}

\DeclareMathOperator{\Edges}{\mathcal{E}}
\newcommand{\av}{\operatorname{av}}
\newcommand{\groot}[1]{(R_{#1}, \chi^{}_{#1})}
\newcommand{\wgrootF}[1]{(R_{#1}, \chi^{}_{#1}, P_{F,#1})} 
\newcommand{\wgroot}[1]{(R_{#1}, \chi^{}_{#1}, P_{#1})} \newcommand{\Fhat}{\widehat{F}}
\newcommand{\Ring}{\mathcal{R}}
\newcommand{\Adm}{\operatorname{Adm}}
\newcommand{\eps}{\varepsilon}
\newcommand{\sublevel}{\mathcal{S}}
\newcommand{\Zhathat}{\widehat{\vphantom{\rule{5pt}{10pt}}\smash{\widehat{Z}}\,}\!}
\newcommand{\Ringt}{\widetilde{\mathcal{R}}}
\newcommand{\Res}{\operatorname{Res}}
\newtheorem{thm}{Theorem}[section]
\newtheorem*{thm*}{Theorem}
\newtheorem{cor}[thm]{Corollary}
\newtheorem{prop}[thm]{Proposition}
\newtheorem{lem}[thm]{Lemma}

\theoremstyle{definition}
\newtheorem{defn}[thm]{Definition}

\newtheorem{exmp}[thm]{Example}

\newtheorem{notn}[thm]{Notation}

\newtheorem{gconvs}[thm]{Grading conventions}

\theoremstyle{remark}
\newtheorem{rem}[thm]{Remark}

\title{Lattice cohomology and  $q$-series invariants of $3$-manifolds}

\begin{document}

\begin{abstract}
An invariant is introduced for negative definite plumbed $3$-manifolds equipped with a spin$^c$-structure. It unifies and extends two theories with rather different origins and structures. One theory is lattice cohomology, motivated by the study of normal surface singularities, known to be isomorphic to the Heegaard Floer homology for certain classes of plumbed 3-manifolds. Another specialization  gives BPS $q$-series which satisfy some remarkable modularity properties and recover ${\rm SU}(2)$ quantum invariants of $3$-manifolds at roots of unity. 
In particular, our work gives rise to a $2$-variable refinement of the $\widehat Z$-invariant. 
\end{abstract}

\author[R. Akhmechet]{Rostislav Akhmechet}
\address{Department of Mathematics, Columbia University, New York NY 10027}
\email{\href{mailto:akhmechet@math.columbia.edu}{akhmechet@math.columbia.edu}}

\author[P. Johnson]{Peter K. Johnson}
\address{Department of Mathematics, Michigan State University, East Lansing, MI 48824}
\email{\href{mailto:john8251@msu.edu}{john8251@msu.edu}}

\author[V. Krushkal]{Vyacheslav Krushkal}
\address{Department of Mathematics, University of Virginia, Charlottesville VA 22904-4137}
\email{\href{mailto:krushkal@virginia.edu}{krushkal@virginia.edu}}

\subjclass[2020]{Primary 57K31; Secondary 57K18, 57K16.}

\maketitle

\section{Introduction}
The development of low-dimensional topology over the last four decades has been greatly influenced by ideas and methods of gauge theory and quantum topology, dating back to the work of Donaldson \cite{Donaldson} and Jones \cite{Jones}. There are many formulations of invariants originating in these theories, but categorically and structurally the two frameworks are quite different: the former is analytic in nature and gives rise to $(3+1)$-dimensional topological quantum field theories, associating to a closed $3$-manifold versions of Floer homology (originally defined in the instanton context in \cite{Floer}). Starting with a quantum group, the latter gives a family of $(2+1)$-dimensional TQFTs, associating to a closed $3$-manifold a collection of numerical  Witten-Reshetikhin-Turaev invariants \cite{Witten, RT} at roots of unity.

Our work builds on two theories that are known to recover, for a certain class of $3$-manifolds, Floer homology and ${\rm SU}(2)$ quantum invariants respectively:
lattice cohomology defined by N\'{e}methi \cite{Nem_Lattice_cohomology} and the $\widehat Z$ invariant of Gukov-Pei-Putrov-Vafa \cite{GPPV}. 
We show that for negative definite plumbed $3$-manifolds, equipped with a spin$^c$-structure, there is a natural construction giving a common refinement of these two theories. As we discuss below, our construction has novel properties that are not satisfied by either lattice cohomology or the $\widehat Z$ invariant. To explain this in more detail, we first summarize the context considered in this paper.

Motivated by the study of normal surface singularities and work of Ozsv\'{a}th-Szab\'{o} \cite{O-S_Floer_Homology_Plumbed}, \cite{Nem_Lattice_cohomology}
introduced lattice cohomology $\Hbb^*(Y, \mathfrak s)$ of negative definite plumbed $3$-manifolds $Y$ with a spin$^c$ structure $\mathfrak s$. For a subclass of negative definite plumbings, $\Hbb^0(Y, \mathfrak s)$ is isomorphic to Heegaard Floer homology $HF^+(-Y, \mathfrak s)$ defined by Ozsv\'{a}th-Szab\'{o} \cite{MR2113019}, as modules over $\Z[U]$. Recent work of Zemke \cite{Zemke} establishes the equivalence between lattice homology and Heegaard Floer homology for general plumbing trees using a completed version of the theories, see Section \ref{sec:lattice cohomology} for a more detailed discussion. The $\Z[U]$-module $\Hbb^0(Y, \mathfrak s)$ was originally defined in \cite{O-S_Floer_Homology_Plumbed}; its generators (as an abelian group) and the action of $U$ are encoded by the {\em graded root}, a certain infinite tree associated to $(Y, \mathfrak s)$, first defined in \cite{Nem_On_the}.

\begin{figure}[t]
\centering
\includegraphics[height=7.7cm]{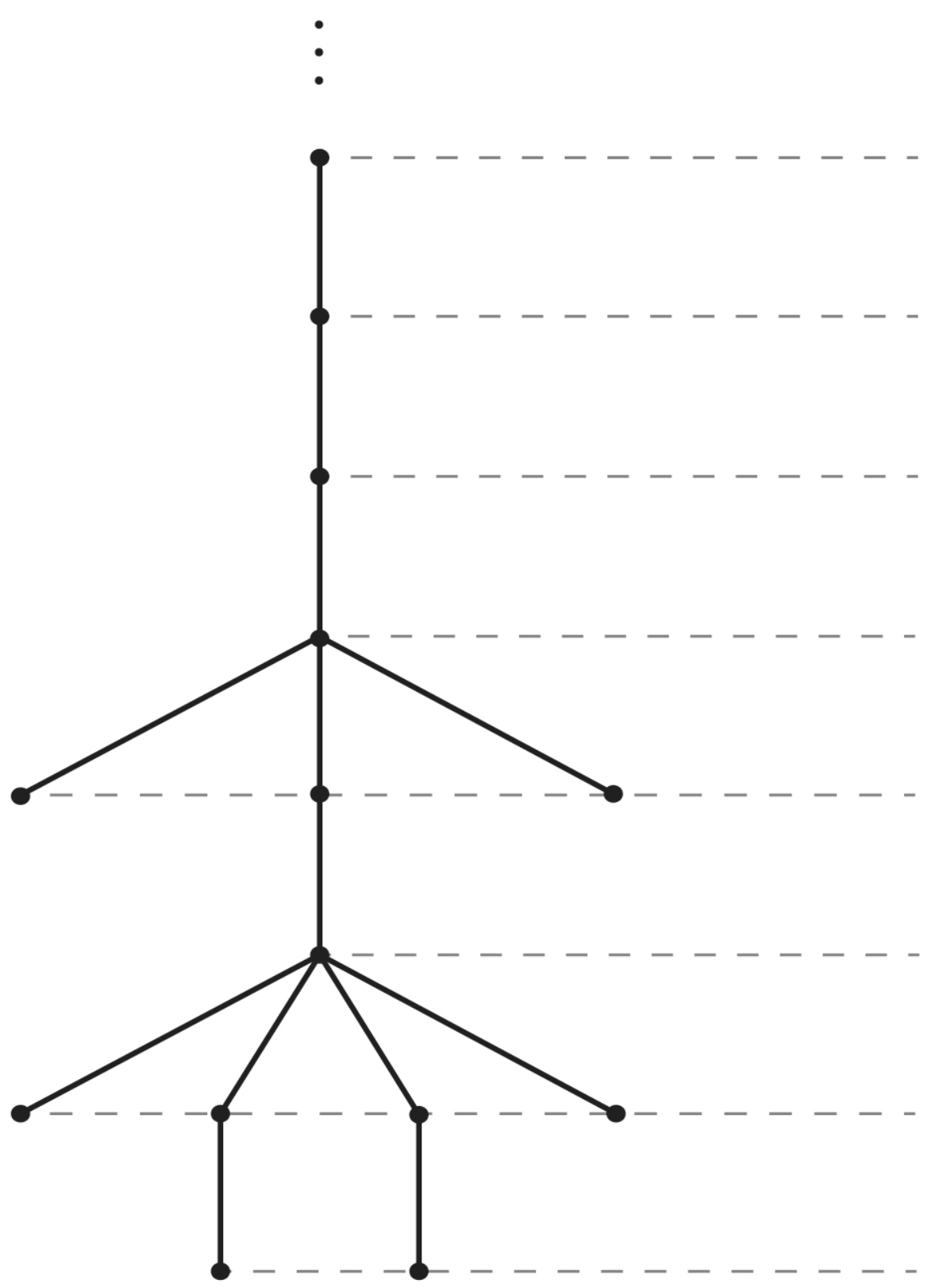}
{\scriptsize
\put(-88,-7){$0$}
\put(-128,-7){$0$}
\put(-164,20){$0$}
\put(-55,20){$0$}
\put(-130,20){$0$}
\put(-86,20){$0$}
\put(-113,57){$0$}
\put(-113,87){$0$}
\put(-55,75){$0$}
}
{\scriptsize
\put(-180,82){$\frac{1}{2}tq^{\frac{13}{2}}$}
\put(-175,110){$\frac{1}{2}(t+t^{-1})q^{\frac{13}{2}}$}
\put(-175,135){$\frac{1}{2}(t+t^{-1})q^{\frac{13}{2}}$}
\put(-200,161){$\frac{1}{2}(t+t^{-1})q^{\frac{13}{2}}- q^{\frac{23}{2}}$}
\put(-200,188){$\frac{1}{2}(t+t^{-1})q^{\frac{13}{2}}- q^{\frac{23}{2}}$}
}
{\small
\put(5,0){$0$}
\put(5,27){$2$}
\put(5,54){$4$}
\put(5,81){$6$}
\put(5,108){$8$}
\put(5,135){$10$}
\put(5,162){$12$}
\put(5,189){$14$}
}
\caption{The weighted graded root for the Brieskorn sphere $\Sigma(2,7,15)$ corresponding to the admissible family $\widehat F$. See grading conventions \ref{grading conv} for an explanation of the numbers in the right column. A more detailed discussion of this example is given in Section \ref{section:conj}.}
\label{Weighted_graded_root} 
\end{figure}
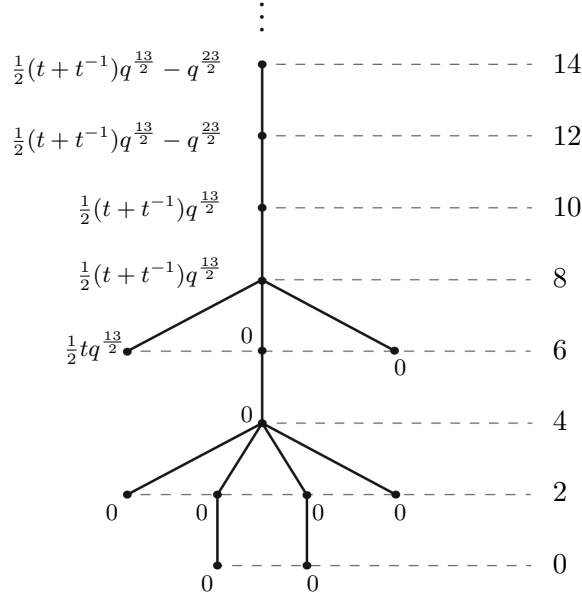

We construct a refinement, an invariant of $(Y,\mathfrak s)$ which takes the form of a graded root labelled by a collection of $2$-variable Laurent polynomials (up to an overall normalization by a fractional power) denoted $P^{}_F$, see Figure \ref{Weighted_graded_root} for an example. As in \cite{Nem_On_the}, the graded root $(R,\chi)$ is defined starting from a negative definite plumbing representing $Y$ and a particular representative of the spin$^c$ structure $\mathfrak s$. The Laurent polynomials $P^{}_F$ labelling the vertices of the graded root in our construction depend on a choice of {\em admissible functions} $F=\{F_n : \Z \to \Ring\}_{n\geq 0}$ where $\Ring$ is a commutative ring, see Definition \ref{def:admissible} for details. 
We work in the setting of $3$-manifolds which are negative definite plumbing trees, as in \cite{Nem_Lattice_cohomology}. 
Our main theorem, proved in Section \ref{sec: The invariant}, is that the result is a topological invariant:

\begin{thm} \label{main thm}
For any admissible family of functions $F$, the weighted graded root $(R, \chi, P^{}_F)$ is an invariant of the $3$-manifold $Y$ equipped with the $\spinc$ structure $\mathfrak s$.
\end{thm}

In Section \ref{sec:two variable series} we show that the sequence of Laurent polynomials $P^n_F$, obtained by summing the labels over the vertices of the graded root $R$ in grading $\chi=n$, stabilizes to a $2$-variable series. Up to an overall normalization, this limit is a Laurent series in $q$ whose coefficients are Laurent polynomials in $t$. Theorem \ref{thm:convergence} shows that it is an invariant of $(Y,\mathfrak s).$

As discussed in Section \ref{sec:weighted graded root}, there is considerable flexibility in the choice of an admissible family of functions $F$. Each such family gives a weighted graded root and a $2$-variable series which are invariants of the $3$-manifold with a spin$^c$ structure. A particular choice, denoted $\widehat F$ in Section \ref{sec:Zhat_a}, gives rise to a $2$-variable series $\Zhathat_{Y,\mathfrak{s}}(q,t)$. To state its properties, we recall the context of the GPPV invariant.

Based on the study of BPS states and certain supersymmetric $6$-dimensional quantum field theories, \cite{GPV, GPPV} formulated a physical definition of homological invariants of $3$-manifolds, denoted $\mathcal{H}^{i,j}_{\text{BPS}} (Y, \mathfrak s)$. 
When the underlying Lie group is ${\rm SU}(2)$, the Euler characterisic of this homology theory is expected to be a $q$-series of the form
\begin{equation} \label{eq: Zhat}
\widehat{Z}_{Y, \mathfrak s}(q) = \sum_{i,j} (-1)^i q^j \, {\rm rk}\,  \mathcal{H}^{i,j}_{\rm{BPS}} (Y, {\mathfrak s}) \; \in \; 2^{-c} q^{\Delta_{\mathfrak s}} \Z[[q]],
\end{equation}
for some $c \in \Z_+$ and $\Delta_{\mathfrak s} \in \Q$ depending on $(Y, {\mathfrak s})$. A mathematically rigorous definition of $\widehat{Z}_{Y,\mathfrak s}(q)$ in general is not yet available. A concrete mathematical formulation for negative definite plumbed $3$-manifolds was given in \cite[Appendix A]{GPPV}; also see Section \ref{sec:Zhat_a} below for a more detailed discussion. An earlier instance of these $q$-series, motivated by the study of WRT invariants and of modular forms, was considered in the case of the Poincar\'{e} homology sphere, and more generally Seifert fibered integer homology spheres with three singular fibers by Lawrence-Zagier in \cite{Lawrence-Zagier}. 
For certain classes of negative definite plumbed $3$-manifolds, the $\widehat Z$ $q$-series are known to satisfy (quantum) modularity properties, cf. \cite{Lawrence-Zagier, MR2757599, 3d_modularity, Quantum_modular_forms}. It is not yet known what kinds of modular forms arise as the $q$-series of other $3$-manifolds including more general negative definite plumbings, and other examples such as Dehn surgeries on hyperbolic knots considered in \cite{GM}.

For Seifert fibered integer homology spheres with three singular fibers, $\widehat Z_{Y,{\mathfrak s}_0}(q)$ is a holomorphic function in the unit disk $|q|<1$, and,  up to a normalization, radial limits to roots of unity give ${\rm SU}(2)$ WRT invariants \cite[Theorem 3]{Lawrence-Zagier}, see also  \cite[Remark 4.5]{GM}. More generally, for rational homology spheres it is conjectured \cite{GPV} that radial limits of a certain linear combination of $\widehat Z_{Y,\mathfrak s}$ over $\spinc$ structures recovers the WRT invariant of $Y$; a precise statement is also given in \cite[Conjecture 3.1]{GM}.

Our next result, established in Sections \ref{sec:two variable series}, \ref{sec:Zhat_a},  relates the $2$-variable series that may be read off from the weighted graded root $(R, \chi, P^{}_{\widehat F})$, as discussed above, to the $\widehat Z$ $q$-series.

\begin{thm} \label{thm: Zhat intro}
The $2$-variable series $\Zhathat_{Y,\mathfrak{s}}(q,t)$ is an invariant of the $3$-manifold $Y$ with a spin$^c$ structure $\mathfrak s$, and its specialization at $t=1$ equals $\widehat{Z}_{Y,\mathfrak s}(q)$.
\end{thm}

The series $\Zhathat_{Y,\mathfrak{s}}(q,t)$ for the Brieskorn sphere $\Sigma(2,7,15)$ in Figure \ref{Weighted_graded_root} is considered in Example  \ref{2 7 15 ex involution0}. It is an interesting question whether there are analogues for the $2$-variable series of the properties of $\widehat Z$ discussed above, in particular the limiting behavior of $\Zhathat(q,t)$ along radial limits of the $q$ variable to roots of unity, as well as modularity of other specializations of $\Zhathat(q,t)$.

Some common features of the $\widehat Z$ invariant with the gauge theory setting were apparent in \cite{GPV, GPPV}; indeed bridging the gap between gauge theoretic and quantum invariants was mentioned as a motivation in \cite{GPPV}. Crucially, the $\widehat Z$ $q$-series depends not just on a $3$-manifold, but also on a spin$^c$ structure. Further, it is shown in \cite{GPPcobordism} that certain numerical gauge-theoretic invariants can be recovered from the $\widehat Z$ series, and a physical discussion of a relation with Heegaard Floer homology is given in \cite{GPV}. Our contribution, as stated in Theorems \ref{main thm}, \ref{thm: Zhat intro} is a new structure that is a common refinement of both perspectives; moreover the weighted graded root $(R, \chi, P^{}_F)$ has new features that are not present in either of them.
Lattice cohomology  and the $\widehat Z$ $q$-series are known to be invariant under conjugation of the $\spinc$ structure, see Section \ref{section:conj} for further details. Corollary \ref{prop: Zhathat conjugation}  states a more subtle transformation of the $2$-variable series $\Zhathat$ under this conjugation. Moreover, Example \ref{prop: Zhathat conjugation} gives a plumbing where conjugate $\spinc$ structures have different weighted graded roots. This example also shows that the Laurent polynomial weights of the graded root carry more information than the limiting series. 

A version of the theory developed here is likely to have an analogue for knot lattice homology of \cite{knot} and the invariant of plumbed knot complements introduced in \cite{GM}. This extension is outside the scope of the present paper; we plan to pursue this in future work.

We conclude by recalling the problem of categorifying WRT invariants of $3$-manifolds, which remains a central open question in quantum topology. The $\widehat Z$ $q$-series provide a very promising approach to this problem. Indeed, as discussed above there is a physical prediction $\mathcal{H}^{i,j}_{\text{BPS}} (Y, \mathfrak s)$ for a homology based on the theory of BPS states \cite{GPV, GPPV}. It is interesting to note that the $2$-variable series $\Zhathat_{Y,\mathfrak{s}}$ constructed in this paper is {\em different} from the expected Poincar\'{e} series of the BPS homology, see Section \ref{sec: Recovering q series}, thus indicating the possibility of a different (or more refined) categorification.

{\bf Acknowledgements.} PJ thanks his advisor, Tom Mark, for his continued support and introducing PJ to lattice cohomology.  VK is grateful to Sergei Gukov for discussions on the GPPV invariant.

RA was supported by NSF RTG grant DMS-1839968, NSF grant DMS-2105467 and the Jefferson Scholars Foundation. PJ was supported by NSF RTG grant DMS-1839968.
VK was supported in part by Simons Foundation fellowship 608604, and NSF Grant DMS-2105467.

\tableofcontents

\newpage

\section{Negative definite plumbed 3-manifolds}

This section summarizes background material and fixes notational conventions  on plumbed $3$-manifolds and $\spinc$ structures; the reader is referred to \cite{Neumann, GompfStipsicz} for more details.

\subsection{Plumbings}  \label{subsec: negative definite plumbings}
\begin{figure}
\centering
\subcaptionbox{A plumbing tree $\Gamma$.
\label{fig:plumbing ex}}[.45\linewidth]
{\begin{tikzpicture}
\node (m) at (0,0) {$\bullet$};
\node (r) at (2,0) {$\bullet$};
\node (l) at (-2,0) {$\bullet$};
\node (b) at (0,-2) {$\bullet$};
\node (s) at (4,0) {$\bullet$};
\node (t) at (6,0) {$\bullet$};

\draw (-2,0) --++ (8,0);
\draw (0,0) --++ (0,-2);

\node[above] at (m) {$-1$};
\node[above] at (r) {$-3$};
\node[above] at (l) {$-2$};
\node[right] at (b) {$-15$};
\node[above] at (s) {$-2$};
\node[above] at (t) {$-2$};

\end{tikzpicture}
}
\subcaptionbox{The framed link $\mathcal{L}(\Gamma)$.
\label{fig:link ex}
}[.5\linewidth]
{\begin{tikzpicture}

\draw (.5,1) arc (90:180:1 and 1);

\draw[preaction={draw, line width=6pt, white}]  (0,0) arc (0:180:1 and 1);

\draw(-2,0) arc (180:360:1 and 1);

\draw(1.5,-1.6) arc (0:90:1 and 1);

\draw[preaction={draw, line width=6pt, white}] (.5,-1) arc (270:360:1 and 1);

\draw[preaction={draw, line width=6pt, white}] (1,0) arc (180:360:1 and 1);

\draw[preaction={draw, line width=6pt, white}] (2.5,0) arc (180:360:1 and 1);

\draw[preaction={draw, line width=6pt, white}] (4,0) arc (180:360:1 and 1);

\draw[preaction={draw, line width=6pt, white}] (-.5,0) arc (180:270:1 and 1);

\draw[preaction={draw, line width=6pt, white}](.5,-.6) arc (90:180:1 and 1);

\draw (-.5,-1.6) arc (180:360:1 and 1);

\draw[preaction={draw, line width=6pt, white}]  (6,0) arc (0:180:1 and 1);

\draw[preaction={draw, line width=6pt, white}]  (4.5,0) arc (0:180:1 and 1);

\draw[preaction={draw, line width=6pt, white}]  (3,0) arc (0:180:1 and 1);

\draw[preaction={draw, line width=6pt, white}] (1.5,0) arc (0:90:1 and 1);

\node[above] at (-1,1) {$-2$};

\node[above] at (2,1) {$-3$};

\node[above] at (3.5,1) {$-2$};

\node[above] at (5,1) {$-2$};

\node[above] at (.5,1) {$-1$};

\node[below] at (.5,-2.6) {$-15$};

\end{tikzpicture}
}
\caption{A negative definite plumbing $\Gamma$ and its associated framed link $\mathcal{L}(\Gamma)$. The $3$-manifold $Y(\Gamma)$ is the Brieskorn sphere $\Sigma(2,7,15)$.}\label{fig:plumbing and link ex}
\end{figure}
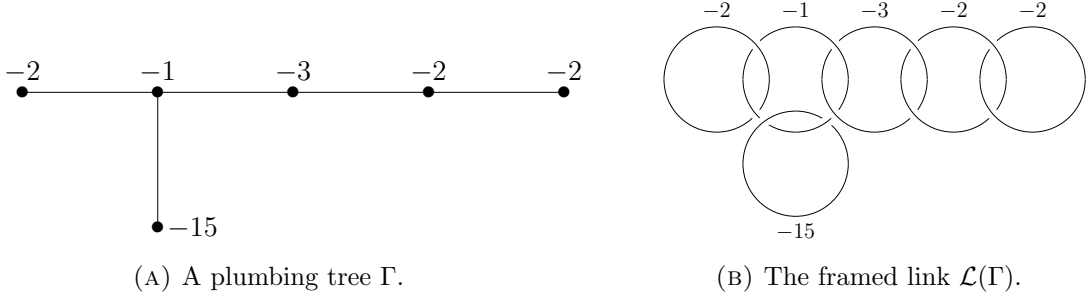

A \emph{plumbing graph} is a finite graph $\Gamma$ equipped with extra data. 
For the purposes of this paper, we restrict to plumbing graphs which are trees equipped with a weight function $m : \V(\Gamma) \to \Z$, where $\V(\Gamma)$ is the set of vertices of $\Gamma$. Let $s= \vert\V(\Gamma)\vert$ be the number of vertices of $\Gamma$ and let $\delta_{v}$ be the degree of $v\in \V(\Gamma)$. We will often implicitly choose an ordering on $\V(\Gamma)$, so that $\V(\Gamma) = \{ v_1, \ldots, v_s\}$, and write quantities associated to $v_i\in \V(\Gamma)$ according to the subscript $i$. For example, $m_i = m(v_i)$, $\delta_i = \delta_{v_i}$, etc. Denote by $m, \delta \in \Z^s$ the weight and degree vectors, respectively:
\begin{equation*}
    m = (m_1, \ldots, m_s), \; \; 
    \delta = (\delta_1, \ldots, \delta_s). 
\end{equation*}
Assign to $\Gamma$ the symmetric $s\times s$ matrix $M= M(\Gamma)$ with entries:
\[
M_{i,j} = 
\begin{cases}
m_{i} & \text{ if } i= j,\\
1 & \text{ if } i\neq j, \text{ and } v_{i} \text{ and } v_{j} \text{ are connected by an edge},\\
0 & {\text{otherwise.} }
\end{cases}
\]
We say $\Gamma$ is \textit{negative definite} if $M$ is negative definite.
From $\Gamma$ we obtain the following manifolds. Consider the framed link $\mathcal{L}(\Gamma)\subset S^3$ given by taking an unknot at each vertex $v$ with framing $m(v)$, and Hopf linking these unknots when the corresponding vertices are adjacent; see Figure \ref{fig:plumbing and link ex} for an example. Let $X=X(\Gamma)$ denote the $4$-manifold obtained by attaching $2$-handles to $D^4$ along $\mathcal{L}(\Gamma)$. Equivalently, $X$ is obtained by plumbing disk bundles over $S^2$ with Euler numbers $m(v)$. Let $Y=Y(\Gamma)$ denote the boundary of $X$, that is the $3$-manifold obtained by Dehn surgery on $\mathcal{L}(\Gamma)$. 

\begin{figure}
\centering
\subcaptionbox{Type (a) move.
\label{fig:type a}}[.5\linewidth]
{\begin{tikzpicture}[scale=.9]
\node(l) at (-2,0) {$\bullet$}; 

\node (m) at (0,0) {$\bullet$};

\node (r) at (2,0) {$\bullet$};

\draw (-2,0) --++ (-1,1);
\draw (-2,0) --++ (-1,.5);
\draw (-2,0) --++ (-1,0);
\draw (-2,0) --++ (-1,-.5);
\draw (-2,0) --++ (-1,-1);

\draw (-2,0) -- (2,0);

\draw (2,0) --++ (1,.5);
\draw (2,0) --++ (1,0);
\draw (2,0) --++ (1,-.5);

\node[anchor = south west] at (l) {$\mathsmaller{m_1 - 1}$};
\node[anchor = south] at (m) {$\mathsmaller{-1}$};
\node[anchor = south east] at (r) {$\mathsmaller{m_2 -1}$};

\node at (-4,0) {$\Gamma' =$};

\begin{scope}[shift = {(0,-3)}]
\node(l) at (-2,0) {$\bullet$};

\node (r) at (2,0) {$\bullet$};

\draw (-2,0) --++ (-1,1);
\draw (-2,0) --++ (-1,.5);
\draw (-2,0) --++ (-1,0);
\draw (-2,0) --++ (-1,-.5);
\draw (-2,0) --++ (-1,-1);

\draw (-2,0) -- (2,0);

\draw (2,0) --++ (1,.5);
\draw (2,0) --++ (1,0);
\draw (2,0) --++ (1,-.5);

\node[anchor = south west] at (l) {$\mathsmaller{m_1}$};

\node[anchor = south east] at (r) {$\mathsmaller{m_2}$};

\node at (-4,0) {$\Gamma  =$};
\end{scope}
\end{tikzpicture}
}
\subcaptionbox{Type (b) move.
\label{fig:type b}
}[.4\linewidth]
{\begin{tikzpicture}[scale=.9]
\node (n) at (0,0) {$\bullet$};
\node (r) at (2,0) {$\bullet$};

\draw (0,0) --++ (-1.5,1);
\draw (0,0) --++ (-1.5,.5);
\draw (0,0) --++ (-1.5,0);
\draw (0,0) --++ (-1.5,-.5);
\draw (0,0) --++ (-1.5,-1);

\draw (0,0) --++ (2,0); 

\node[anchor = south west] at (n) {$\mathsmaller{m_1 -1}$};
\node[anchor = south] at (r) {$\mathsmaller{-1}$};
\node at (-2.4,0) {$\Gamma' =$};

\begin{scope}[shift= {(0,-3)}]
\node (n) at (0,0) {$\bullet$};

\draw (0,0) --++ (-1.5,1);
\draw (0,0) --++ (-1.5,.5);
\draw (0,0) --++ (-1.5,0);
\draw (0,0) --++ (-1.5,-.5);
\draw (0,0) --++ (-1.5,-1);

\node[anchor = south west] at (n) {$\mathsmaller{m_1}$};
\node at (-2.4,0) {$\Gamma  =$};
\end{scope}
\end{tikzpicture}
}
\caption{Type (a)  and (b) Neumann moves relating negative definite plumbing trees for homeomorphic $3$-manifolds.  }\label{fig:Neumann moves}
\end{figure}
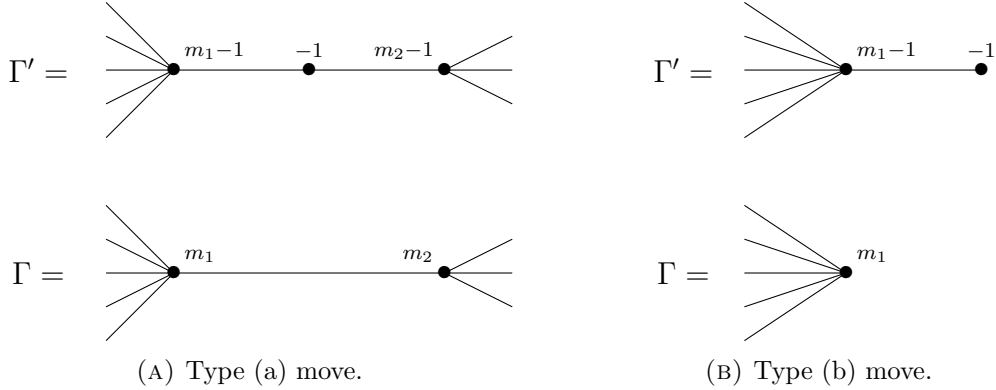

A \emph{negative definite plumbed 3-manifold} is a 3-manifold that is homeomorphic to $Y(\Gamma)$ for some negative definite plumbing tree $\Gamma$. In \cite{Neumann}, the relationship between different representations of a given negative definite plumbed 3-manifold is studied. In particular, from the results in \cite{Neumann}, one can deduce the following theorem which is used in \cite[Proposition 3.4.2]{Nem_Lattice_cohomology} to prove the topological invariance of lattice cohomology:

\begin{thm}[\cite{Neumann}] \label{Neumann thm}
Two negative definite plumbing trees $\Gamma$ and $\Gamma'$ represent the same 3-manifold if and only if they can be related by a finite sequence of type (a) and (b) Neumann moves shown in Figure \ref{fig:Neumann moves}. 
\end{thm}

\begin{notn}
\label{notation type a and b}
For future reference, we establish notation associated to type (a) and (b) Neumann moves.

\begin{itemize}
    \item We use primes to distinguish quantities associated with $\Gamma'$ from those associated with $\Gamma$. For example, $\delta'$ denotes the degree vector for $\Gamma'$, whereas $\delta$ denotes the degree vector for $\Gamma$.
    \item For a type (a) move, we will always order vertices so that the $-1$ weighted vertex in $\Gamma'$ which is blown down is labeled by $v_{0}'$, and the two adjacent vertices with weights $m_{1}-1$ and $m_{2}-1$ are labeled by $v_{1}'$ and $v_{2}'$ respectively. In the $\Gamma$ graph of a type (a) move, there is no vertex $v_{0}$ and the two vertices with weights $m_{1}$ and $m_{2}$ are labeled by $v_{1}$ and $v_{2}$ respectively.
    
    \item For a type (b) move, the $-1$ weighted vertex on $\Gamma'$ is labeled by $v_{0}'$ and its adjacent vertex is labeled by $v_{1}'$. In the $\Gamma$ graph, there is no $v_{0}$ vertex and the vertex with weight $m_{1}$ is labeled by $v_{1}$.
\end{itemize}
\end{notn}

\subsection{Identification of \texorpdfstring{$\spinc$}{} structures}\label{sec:identification of spinc structures}
$\Spinc$ structures are important ingredients to both lattice cohomology and the $\Zhat$-invariant. The two theories use different identifications of $\spinc$ structures in terms of plumbing data. We recall a translation between the two identifications, following \cite[Section 4.2]{GM}.

To begin, we describe the relationships between various (co)homology groups of $X$ and $Y$. First, note that $\Gamma$ gives a convenient choice of basis for $H_2(X;\Z)$ in the following way. For $v\in \V(\Gamma)$, let $[v]\in H_2(X;\Z)$ be the class of the $2$-sphere obtained by capping off the core of the $2$-handle corresponding to $v$. Then $H_2(X;\Z) \cong \Z^s$, with a basis given by $\{[v_{1}], \ldots, [v_{s}]\}$. With respect to this basis, the intersection form on $H_2(X;\Z)$
is the bilinear form associated with $M$, $\langle x,y\rangle = x^t M y$. We will also write
\[
\langle-,-\rangle : \Z^s \times \Z^s \to \Z
\]
to denote this bilinear form when $H_2(X;\Z)$ is identified with  $\Z^s$ as above.

\begin{rem}
In some of the lattice cohomology literature the intersection form is denoted by $(-,-)$. However, in \cite{GM} the intersection form is denoted using angled brackets $\langle-,-\rangle$, as we do above, and $(-,-)$ instead refers to the usual dot product. To minimize confusion, we will use $\cdot$ for the dot product.
\end{rem}

Since $X$ is a 2-handlebody with no 1-handles, $H_1(X;\Z) = 0$, and we can identify $H^2(X;\Z)$ with $\Hom_{\Z}(H_2(X;\Z),\Z)$. Furthermore, using the above basis of $H_2(X;\Z)$, we have a distinguished isomorphism $\Hom_{\Z}(H_2(X;\Z),\Z) \cong \Z^s$, so that a vector $k\in \Hom_{\Z}(H_2(X;\Z),\Z)$ is represented as $k = (k([v_1]), \ldots, k([v_s])) \in \Z^s$. Combining these two identifications, we get an identification of $H^2(X;\Z)$ with $\Z^s$ such that for $k\in H^2(X,\Z)$ and $x\in H_2(X;\Z)$, we have $k(x) = k \cdot x$. The identifications described above are used throughout the paper.

\begin{defn}
An element $k\in H^2(X;\Z)$ is called \textit{characteristic} if $k\cdot x + \langle x, x\rangle \equiv 0 \bmod 2$ for all $x\in H_{2}(X;\Z)$. We denote the set of characteristic vectors of $X$ by $\Char(X)$. 
\end{defn}
In terms of our identification of $H^2(X;\Z)$ with $\Z^s$, it follows that: 
\begin{align*}
    \Char(X) = m + 2\Z^s.
\end{align*}
It is a standard fact in 4-manifold topology that for simply connected $X$ the map which takes a $\spinc$ structure on $X$ to the first Chern class of its determinant line bundle is a bijection from $\spinc(X)$ to $\Char(X)$, cf. \cite[Proposition 2.4.16]{GompfStipsicz}. Moreover, by restricting $\spinc$ structures to the boundary 3-manifold $Y$, we get the following identification:
\begin{equation}
\label{eq:LC spinc structures}
    \spinc(Y) \cong \frac{m+2\Z^s}{2M\Z^s}. 
\end{equation}
The right side of the above identification is to be interpreted as the set $m+2\Z^s$ up to the equivalence relation defined by $k\sim k'$ if $k-k'\in 2M\Z^s$. If $k\in m +2\Z^s$, we let $[k]\in \spinc(Y)$ denote the equivalence class containing $k$. 

The identification of $\spinc$ structures just described is the one used in the lattice cohomology and Heegaard Floer homology literature (see for example \cite[Section 2.2.2]{Nem_Lattice_cohomology} and \cite[Section 1]{O-S_Floer_Homology_Plumbed}). The identification used in the $\Zhat$ literature is given as follows:
\begin{equation}
    \label{eq:GM spinc structures}
    \spinc(Y) \cong \frac{\delta+2\Z^s}{2M\Z^s}. 
\end{equation}
Note that in this identification, unlike in \eqref{eq:LC spinc structures}, the set $\delta +2\Z^s$ is not necessarily equal to the set of characteristic vectors. In fact, $\delta +2\Z^s = \Char(X)$ if and only if $\delta\equiv m \bmod{2}$.

As in identification \eqref{eq:LC spinc structures}, for $a\in \delta + 2\Z^s$, we let $[a]\in \spinc(Y)$ denote the equivalence class containing $a$. To avoid confusion between the two $\spinc$ identifications, throughout the paper we will use the letter $k$ for vectors in $m + 2\Z^s$ and the letter $a$ for vectors in $\delta + 2\Z^s$.

The lattice cohomology and $\Zhat$ identifications of $\spinc$ structures are related to each other in the following way. For a fixed plumbing graph $\Gamma$, there is a bijection:
\begin{equation}
    \label{eq:translating between spinc identifications}
      \frac{m+2\Z^s}{2M\Z^s}   \overset{\sim}{\longrightarrow} \frac{\delta+2\Z^s}{2M\Z^s}, \hspace{ 1cm}
      [k]  \longmapsto [k-(m+\delta)]. 
\end{equation}

\begin{rem}
If we let $u = (1, 1,\ldots, 1)\in \Z^s$, then $Mu = m+\delta$ and we can alternatively express the above bijection as $[k] \mapsto [k-Mu]$.
\end{rem}

The above bijection is natural with respect to type (a) and type (b) Neumann moves in the sense that if $\Gamma$ and $\Gamma'$ are two plumbing graphs related by either one of these two moves, we get the following commutative diagram:

\begin{center}
\begin{tikzpicture}[baseline= (a).base]
\node[scale=1.2] (a) at (0,0){
\begin{tikzcd}
\large
     \frac{m+2\Z^s}{2M\Z^s} 
\ar[r] \ar{d}{\alpha} & \frac{\delta+2\Z^s}{2M\Z^s} \ar{d}{\beta}\\
        \frac{m'+2\Z^{s+1}}{2M'\Z^{s+1}} \ar[r] & \frac{\delta'+2\Z^{s+1}}{2M'\Z^{s+1}}
\end{tikzcd}
};
\end{tikzpicture}
\end{center}

Here $\alpha, \beta$  are also bijections, which at the level of representatives are defined as follows. 

Type (a) move:
\begin{equation}
    \label{eq:k' for type a}
       \alpha(k) =k':= (0, k) + (1, -1, -1, 0, \ldots, 0), \hspace{.5cm} 
       \beta(a)=a':=(0, a). 
\end{equation}

Type (b) move: 
\begin{equation}
    \label{eq:k' for type b}
      \alpha(k) =k':=(0, k) + (-1, 1, 0,\ldots, 0), \hspace{.5cm}
      \beta(a)=a':=(0, a) + (-1, 1, 0,\ldots, 0).
\end{equation}

\section{Lattice cohomology}
\label{sec:lattice cohomology}
Lattice cohomology was introduced by Némethi in \cite{Nem_Lattice_cohomology}, building on earlier work of Ozsváth-Szabó in \cite{O-S_Floer_Homology_Plumbed}. It is a theory which assigns to a given plumbing graph $\Gamma$ and $\spinc$ structure $[k]$, a $\Z[U]$-module:
\begin{align*}
   \Hbb^*(\Gamma, [k]) = \bigoplus_{i=0}^{\infty}\Hbb^i(\Gamma, [k]) 
\end{align*}
Each $\Hbb^i(\Gamma, [k])$ is a $(2\Z)$-graded $\Z[U]$-module. Hence, $\Hbb^*$ is bigraded; it carries a homological grading given by the superscript $i$ as well as an internal $(2\Z)$-grading. 

It was shown in \cite{Nem_Lattice_cohomology} that for negative definite plumbings $\Hbb^*$ is invariant under Neumann moves and therefore is a topological invariant. Furthermore, extending results from \cite{O-S_Floer_Homology_Plumbed}, it was shown that for a subset of negative definite plumbed 3-manifolds, namely \textit{almost rational} plumbings, there exists an isomorphism between lattice cohomology and Heegaard Floer homology:

\begin{thm}[{\cite[Theorems 4.3.3 and 5.2.2]  {Nem_Lattice_cohomology}}]\label{thm: LC-HF iso}
If $\Gamma$ is almost rational, then as graded $\Z[U]$-modules,
\begin{align}
\begin{aligned}
 \label{eq:LC-HF iso}
    \Hbb^i(\Gamma, [k])\left[-\max\limits_{k'\in [k]}\frac{(k')^2+s}{4}\right]&\cong 
    \begin{cases}
    HF^+(-Y(\Gamma), [k])&\text{if }i =0\\
    0&\text{otherwise}.
    \end{cases}
\end{aligned}
\end{align}
Here $-Y$ is $Y$ with the opposite orientation, $(k')^2:=(k')^tM^{-1}k'$, and the left side of the equation is $\Hbb^i(\Gamma, [k])$ with its internal grading shifted up by $-\max\limits_{k'\in [k]}\frac{(k')^2+s}{4}$. 
\end{thm}

\begin{rem}\label{rem: k^2 meaning}
The quantity $(k')^2$ has the following geometric meaning: it is the square of the first Chern class of the $\spinc$ structure on $X(\Gamma)$ corresponding to the characteristic vector $k'$. Even if a vector $x$ is not characteristic, we will still define $x^2 := x^tM^{-1}x\in \Q$.
\end{rem}

\begin{rem}\label{rem: d inv}
The minimal internal grading of $\Hbb^0(\Gamma, [k])$ is always equal to zero. Hence, after the grading shift, the minimal grading of the left side of equation \eqref{eq:LC-HF iso} is equal to $-\max\limits_{k'\in [k]}\frac{(k')^2+s}{4}$. In particular, Theorem \ref{thm: LC-HF iso} implies the $d$-invariant of $-Y$ at the $\spinc$ structure $[k]$ is equal to $-\max\limits_{k'\in [k]}\frac{(k')^2+s}{4}$, or, equivalently, the $d$-invariant of $Y$ at $[k]$ is $\max\limits_{k'\in [k]}\frac{(k')^2+s}{4}$. 
\end{rem}

In light of Theorem \ref{thm: LC-HF iso}, we focus our attention on $\Hbb^0$ rather than recalling the full definition of lattice cohomology. However, it is worth noting that there do exist negative definite plumbings with non-trivial $\Hbb^{i}$, $i\geq 1$ (see for example \cite[Example 4.4.1]{Nem_Lattice_cohomology}).
Moreover, there is a related theory $\Hbb_*(\Gamma, \mathfrak{s})$ of lattice \emph{homology}  \cite{Nem_Lattice_cohomology, OSS-spectral-sequence}, which can be defined for not necessarily negative definite plumbing trees if one works with a completed version of the construction over the power series ring $\mathbb{F}[[U]]$, where $\mathbb{F}=\Z/2\Z$. The equivalence of $\Hbb_*(\Gamma, \mathfrak{s})$ and (a completed version of) $HF^-(Y(\Gamma),\mathfrak{s})$ was recently established in \cite{Zemke}, building on work of \cite{OSS-spectral-sequence, Zemke_bordered}.

\subsection{Graded roots}
The $(2\Z)$-graded $\Z[U]$-module $\Hbb^0(\Gamma, [k])$ has the nice feature that it can be encoded by a graph, called a \textit{graded root}, proven in \cite[Proposition 4.6]{Nem_On_the} to itself be a topological invariant. We now recall the notion of a graded root $(R, \chi)$ as an abstract object and show how to associate to it a $(2\Z)$-graded $\Z[U]$-module $\Hbb(R,\chi)$. In the next subsection, we show how to obtain a graded root $\groot{[k]}$ from a pair $(\Gamma, [k])$ and we define $\Hbb^0(\Gamma, [k]) = \Hbb\groot{[k]}$. For complete details, see \cite[Section 3]{Nem_On_the}.

\begin{figure}
    \centering
    \begin{tikzpicture}[scale=0.5]
\node at (0,0) {$\bullet$};
\node at (2,0) {$\bullet$};
\node at (1,1.5) {$\bullet$};
\node at (3,1.5) {$\bullet$};
\node at (2,3) {$\bullet$};
\node at (2,4.5) {$\bullet$};
\node at (2,6) {$\bullet$};

\draw (0,0)--(1,1.5);
\draw (2,0)--(1,1.5);
\draw (1,1.5)--(2,3);
\draw (3,1.5)--(2,3);
\draw (2,3)--(2,4.5)--(2,6);

\node at (4.5,0) {$\chi = 5$};
\node at (4.5,1.5) {$\chi= 6$};
\node at (4.5,3) {$\chi = 7$};
\node at (4.5,4.5) {$\chi = 8$};
\node at (4.5,6) {$\chi =9$};
\node at (4.5, 7) {$\vdots$};
\end{tikzpicture}
    \caption{An example of a graded root.} 
    \label{fig:graded root example}
\end{figure}
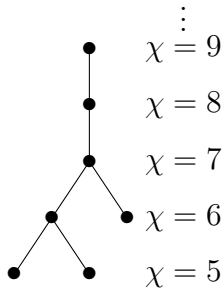

\begin{defn}
\label{def:graded root} A \emph{graded root} is an infinite tree $R$, with vertices and edges denoted $\V$ and $\Edges$ respectively, together with a grading function $\chi : \V \to \Z$ satisfying the properties listed below. We write an edge with endpoints $u$ and $v$ as $[u,v]\in \Edges$.
\begin{enumerate}
    \item $\chi(u) - \chi(v) = \pm 1$ for any $[u,v] \in \Edges$. 
    \item $\chi(u) > \min \{ \chi(v), \chi(w) \}$ for any $[u,v], [u,w] \in \Edges$ with $v\neq w$. 
    \item $\chi$ is bounded below, each preimage $\chi^{-1}(n)$ is finite, and $\vert \chi^{-1}(n) \vert =1$ for sufficiently large $n$. 
\end{enumerate}
An isomorphism of graded roots is an isomorphism of the underlying graphs that respects the grading. For $r\in \Z$, let $(R,\chi)[r] = (R,\chi[r])$ denote the graded root with the same underlying tree and the grading shifted up by $r$, so that $\chi[r](v) = \chi(v) + r$. 
 
\end{defn}
We visualize a graded root by embedding it into the plane such that vertices of the same grading are placed at the same horizontal level, see Figure \ref{fig:graded root example} for an example. 

The $(2\Z)$-graded $\Z[U]$-module $\Hbb(R,\chi)$ associated to a graded root $(R, \chi)$ is defined as follows:
\begin{itemize}
    \item To each $v\in\V$, we associate a copy of $\Z$, which we denote $\Z_v$. By an abuse of notation, we let $v$ also denote a distinguished generator of $\Z_{v}$.
    \item As a graded $\Z$-module, we let $\Hbb(R,\chi):=\bigoplus\limits_{v\in \V}\Z_{v}$, where $\Z_{v}$ has grading equal to $2\chi(v)$.
    \item For each generator $v$, we let $Uv = v_{1}+\cdots + v_{n}$ where 
    \begin{align*}
        \{v_{1},\ldots, v_{n}\} = \{w\in \V\ |\ \chi(w) = \chi(v)-1\text{ and } w\text{ is connected to }v\text{ by an edge}\}.
    \end{align*}
    In particular, if the above set is empty, then $Uv =0$.
    \item We then extend the $U$-action by $\Z$-linearity. Note $U$ decreases grading by $2$.
    
\end{itemize}

\subsection{Graded roots for negative definite plumbings}
\label{sec:graded roots for neg def plumbings}
We now recall from \cite[Section 4]{Nem_On_the} the graded root $\groot{k}$ associated to a negative definite plumbing $\Gamma$ and a $\spinc$ representative $k\in m + 2\Z^s$. We then show how to obtain a graded root $\groot{[k]}$ corresponding to $(\Gamma, [k])$, independent of the choice of $\spinc$ representative $k$.

Define a function $\chi^{}_k : \Z^s \to \Z$ by
\begin{equation}
\label{eq:chi}
    \chi^{}_k(x) = - (k\cdot x + \langle x,x\rangle)/2.
\end{equation}
Note that $\chi^{}_k(x)\in \Z$ since $k$ is characteristic.

Consider the standard cubical complex structure on $\R^s$, with $0$-dimensional cells located at the lattice points $\Z^s \subset \R^s$. We can extend $\chi^{}_k$ to a function on closed cells $\square$ (of any dimension), by defining
\begin{align*}
    \chi^{}_{k}(\square) = \max\{\chi^{}_{k}(v) \mid v\text{ is a 0-cell of }\square\}
\end{align*}

Let $S_j \subset \R^s$ denote the subcomplex consisting of cells $\square$ such that $\chi^{}_k(\square) \leq j$. We call $S_j$ a \emph{sublevel set}. Note that each $S_{j}$ is compact since the intersection form $\langle-,-\rangle$ is assumed to be negative definite. More precisely, if one considers $\chi^{}_{k}$ as a function from $\R^s\to \R$, then it is bounded below and its level sets are $(s-1)$-dimensional ellipsoids. 

Write each sublevel set as a disjoint union over its connected components,
\begin{align*}
    S_j = C_{j,1} \sqcup \cdots \sqcup C_{j,n_j}
\end{align*}
The vertices of $R_k$ consist of connected components among all the $S_j$,
\[
\V := \{C_{j,\ell} \mid j\in \Z, 1\leq \ell \leq n_j\},
\]
and the grading is given by $\chi^{}_k(C_{j,\ell}) = j$, where as in \cite{Nem_On_the} we use $\chi_{k}$ to denote both a function on closed cells of our cellular decomposition as well as a grading on the vertices of $R_k$. Edges of $R_k$ correspond to inclusions of connected components: there is an edge connecting $C_{j,\ell}$ and $C_{j+1, \ell'}$ if $C_{j,\ell} \subseteq C_{j+1,\ell'}$. By \cite[Proposition 4.3]{Nem_On_the}, $\groot{k}$ is a graded root. 

Let us now recall from \cite{Nem_On_the} how $\groot{k}$ depends on the choice of a $\spinc$ representative. Let $k\in m+2\Z^s$ and let $k'=k + 2My$ be another representative for $[k]\in \spinc(Y)$. One readily checks that: 
\begin{equation}
\label{eq:chi for different reps}
    \chi^{}_{k'}(x) = \chi^{}_{k}(x+y) - \chi^{}_{k}(y)
\end{equation}
for all $x\in \Z^s$. As stated in \cite[Proposition 4.4]{Nem_On_the}, there is an isomorphism of graded roots: 
\begin{equation}\label{eq:GR for different reps}
\groot{k'} \cong \groot{k}[-\chi^{}_k(y)],
\end{equation}
given by applying the translation $x\mapsto  x+y$, to each connected component $C$ in each sublevel set of $\chi^{}_{k'}$. Since the collection of graded roots $\{\groot{k'}\}_{k'\in [k]}$ are all isomorphic up to an overall grading shift, we normalize gradings in the following way to obtain a graded root independent of the choice of $\spinc$ representative. This is the same normalization as \cite[Section 4.5.1]{Nem_On_the}. 

\begin{defn}
\label{def:normalizating graded root}
Let $[k]\in \spinc$ and define $\groot{[k]}$ by taking any representative $k'\in [k]$ and shifting the $\chi^{}_{k'}$-grading on $(R_{k'}, \chi^{}_{k'})$ so that its minimal grading is zero.
\end{defn}

\begin{gconvs}\label{grading conv}
When drawing the graded root $\groot{[k]}$ associated to a plumbing $\Gamma$ and $\spinc$ structure $[k]$, the numbers we list in the vertical column to the right are the gradings of the corresponding generators of $\Hbb^0(\Gamma, [k])\left[-\max\limits_{k'\in [k]}\frac{(k')^2+s}{4}\right]$. The reason we do this is so that when the isomorphism (\ref{eq:LC-HF iso}) holds, the gradings one sees are the $HF^+$ gradings. See for example Figure \ref{Weighted_graded_root}, where $d(-\Sigma(2,7,15))=0$. 
\end{gconvs}

\section{Admissible functions and weighted graded roots}
\label{sec:weighted graded root}

This section illustrates the main construction of the paper in a preliminary context.
Let $\Gamma$ be a negative definite plumbing and $k\in m + 2\Z^s$ a $\spinc$ representative. Given a function 
\begin{equation}
\label{eq:F Gamma k}
F_{\Gamma,k} : \Z^s \to \Ring
\end{equation}
valued in some ring $\Ring$, each vertex $v$ in the graded root $\groot{k}$ can be given a {\em weight} by taking the sum of $F_{\Gamma,k}$ over lattice points in the connected component representing $v$. Precisely, for a connected component $C$ in some sublevel set of $\chi^{}_k$, let $L(C) = C \cap \Z^s$ denote its lattice points, and define its weight to be
\begin{equation}
\label{eq:sum over connected component}
F_{\Gamma,k}(C) := \sum_{x\in L(C)} F_{\Gamma,k}(x).
\end{equation}

Note that $L(C)$ is finite since the sublevel sets are compact. To obtain an invariant of $3$-manifolds, the weights $F_{\Gamma,k}(C)$ should be invariant under the Neumann moves in Theorem \ref{Neumann thm}; this imposes significant constraints on the functions \eqref{eq:F Gamma k}.
In this section we explain a way to obtain $F_{\Gamma,k}$ satisfying these constraints from an \emph{admissible} family of functions $F = \{F_n\}_{n\geq 0}$. Theorem \ref{thm:sum of weights is invt} shows that the graded root $\groot{k}$ with these weights is an invariant of $(Y(\Gamma), [k])$.
This result follows from the more general Theorem \ref{thm:invariance}.

\begin{defn} \label{def:admissible}  Fix a commutative ring $\Ring$.
A family of functions $F = \{F_n : \Z \to \Ring\}_{n\geq 0}$ is \emph{admissible} if 
\begin{enumerate}[label= (A\arabic*)]
    \item\label{item:A1} $F_2(0) = 1$ and $F_2(r) = 0$ for all $r\neq 0$. 
    \item \label{item:A2} For all $n\geq 1$ and $r\in \Z$, 
    \[ F_n(r+1) - F_n(r-1) = F_{n-1}(r).\]
\end{enumerate}
\end{defn}

 For an admissible family $F = \{F_n\}_{n\geq 0}$, define $F_{\Gamma,k} :\Z^s \to \Ring$ by
\begin{equation}
\label{eq:F_Gamma,k}
    F_{\Gamma,k}(x) = \prod_{v\in \V(\Gamma)} F_{\delta_v}\left( (2Mx + k - Mu)_v \right),
\end{equation}
where $u=(1,1,\ldots, 1)\in \Z^s$ and $(-)_v$ denotes the component corresponding to $v$.

\begin{rem} \label{Gamma remark}
We stress that $F_n$ and $F_{\Gamma,k}$ are only set-theoretic functions rather than homomorphisms. The definition \eqref{eq:F_Gamma,k} of $F_{\Gamma,k}(x)$ is motivated by the $\widehat Z$-invariant, see Section \ref{sec:Zhat_a} and in particular Proposition \ref{prop:Zhat_k in terms of F}.
\end{rem}

Note that if $k'=k+2My$ is another representative for $[k]\in \spinc(Y(\Gamma))$, then 
\begin{equation}
\label{eq:F different representatives}
    F_{\Gamma,k'}(x) = F_{\Gamma,k}(x+y),
\end{equation}
so the weights in equation \eqref{eq:sum over connected component} are compatible with the isomorphism $\groot{k'} \cong \groot{k}[-\chi^{}_k(y)]$ from Section \ref{sec:graded roots for neg def plumbings}. Denote by $(R_{[k]}, \chi^{}_{[k]}, F_{[k]})$ the graded root $\groot{[k]}$ equipped with these weights.

\begin{thm}
\label{thm:sum of weights is invt}
For any admissible family of functions $F$, the weighted graded root $(R_{[k]}, \chi^{}_{[k]}, F_{[k]})$ is an invariant of the $3$-manifold $Y(\Gamma)$ endowed with the $\spinc$ structure $[k]$. 
\end{thm}
\noindent The proof of this result follows from Theorem \ref{thm:invariance} upon specializing $q=t=1$. 

We pause to make some remarks about admissible families of functions. Explicit $\Z$- and $\Z[\frac{1}{2}]$-valued examples motivated by the $\Zhat$ invariant are given in Definitions \ref{def:Fhat} and \ref{def:Fhatpm}. Note that not only $F_2$, but also $F_0$ and $F_1$ are uniquely determined by conditions \ref{item:A1} and \ref{item:A2}: 
\begin{equation}
\label{eq:F_1 and F_0}
F_1(r) =
\begin{cases}
1 & \text{ if } r=-1, \\
-1 & \text{ if } r=1, \\
0 & \text{ otherwise.} 
\end{cases}
\hskip3em 
F_0(r) =
\begin{cases}
1 & \text{ if } r=\pm 2, \\
-2 & \text{ if } r=0, \\
0 & \text{ otherwise.} 
\end{cases}
\end{equation}

We note that each factor $F_{\delta_v}((2Mx + k - Mu)_v)$ in equation \eqref{eq:F_Gamma,k} depends on $x$ only via $(Mx)_v$, the coordinate corresponding to $v$, so level sets of $F_{\delta_v}((2Mx + k - Mu)_v)$ are hyperplanes. 
Therefore $F_{\Gamma,k}$ is supported on finitely many hyperplanes when some $F_{\delta_v}$ has finite support. If every vertex of $\Gamma$ has degree at most $j$, then only the functions $F_i$ for $i\leq j$ appear in the definition of $F_{\Gamma,k} : \Z^s\to \Ring$ from equation \eqref{eq:F_Gamma,k}. In particular, if every vertex of $\Gamma$ has degree at most two (plumbings of this form are lens spaces), then  \eqref{eq:F_1 and F_0} and property \ref{item:A1} imply that $F_{\Gamma,k} $ has finite support. 
In general, $F_{\Gamma,k}$ need not have finite support if $\Gamma$ contains a vertex of degree at least three, see for example the admissible family $\h{F}$ from Definition \ref{def:Fhat} and the related admissible families $\h{F}^{\pm}$ of Definition \ref{def:Fhatpm}.

We can characterize admissible families in the following way. Let $\Adm(\Ring)$ denote the set of all admissible $\Ring$-valued families of functions, and let $(\Ring \times \Ring)^{\mathbb{N}}$ be the set of all sequences with entries in $\Ring \times \Ring$. 

\begin{prop}\label{prop: adm identification}
There is a bijection $\Adm(\Ring) \cong (\Ring \times \Ring)^{\mathbb{N}}$.
\end{prop}

\begin{proof}
We show that $\Psi: \Adm(\Ring)\to (\Ring \times \Ring)^{\mathbb{N}}$ defined by $\Psi(F) = (F_{n+2}(0), F_{n+2}(1))_{n\geq 1}$ is a bijection. Suppose $F = \{F_{n}\}_{n\geq 0}$ is an admissible family. Fix $n\geq 1$. By applying the recursive relation \ref{item:A2} inductively, we see that $F_n(r)$ is determined by $F_{n-1}, F_{n}(0)$ if $r$ is even and $F_{n-1},F_{n}(1)$ if $r$ is odd, so $\Psi$ is injective. Likewise, given $(a_{n}, b_{n})_{n\geq 1}\in (\Ring \times \Ring)^{\mathbb{N}}$, we set $(F_{n+2}(0),F_{n+2}(1))=(a_n, b_n)$ and use \ref{item:A2} to inductively construct an admissible family $F$ such that $\Psi(F) = (a_{n}, b_{n})_{n\geq 1}$.
\end{proof}

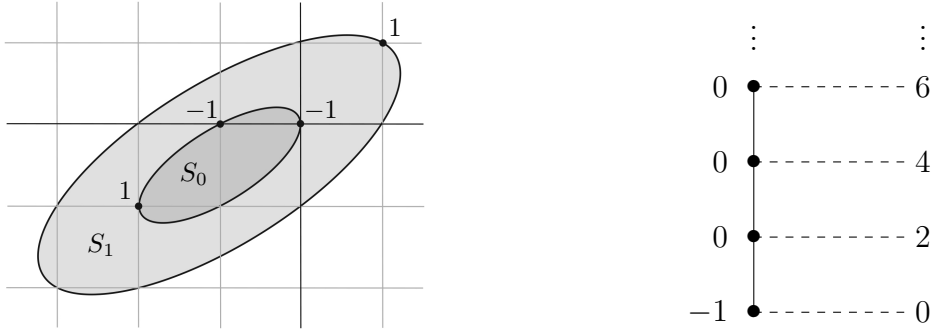
\begin{figure}
    \centering
{\includegraphics[height=4.5cm]{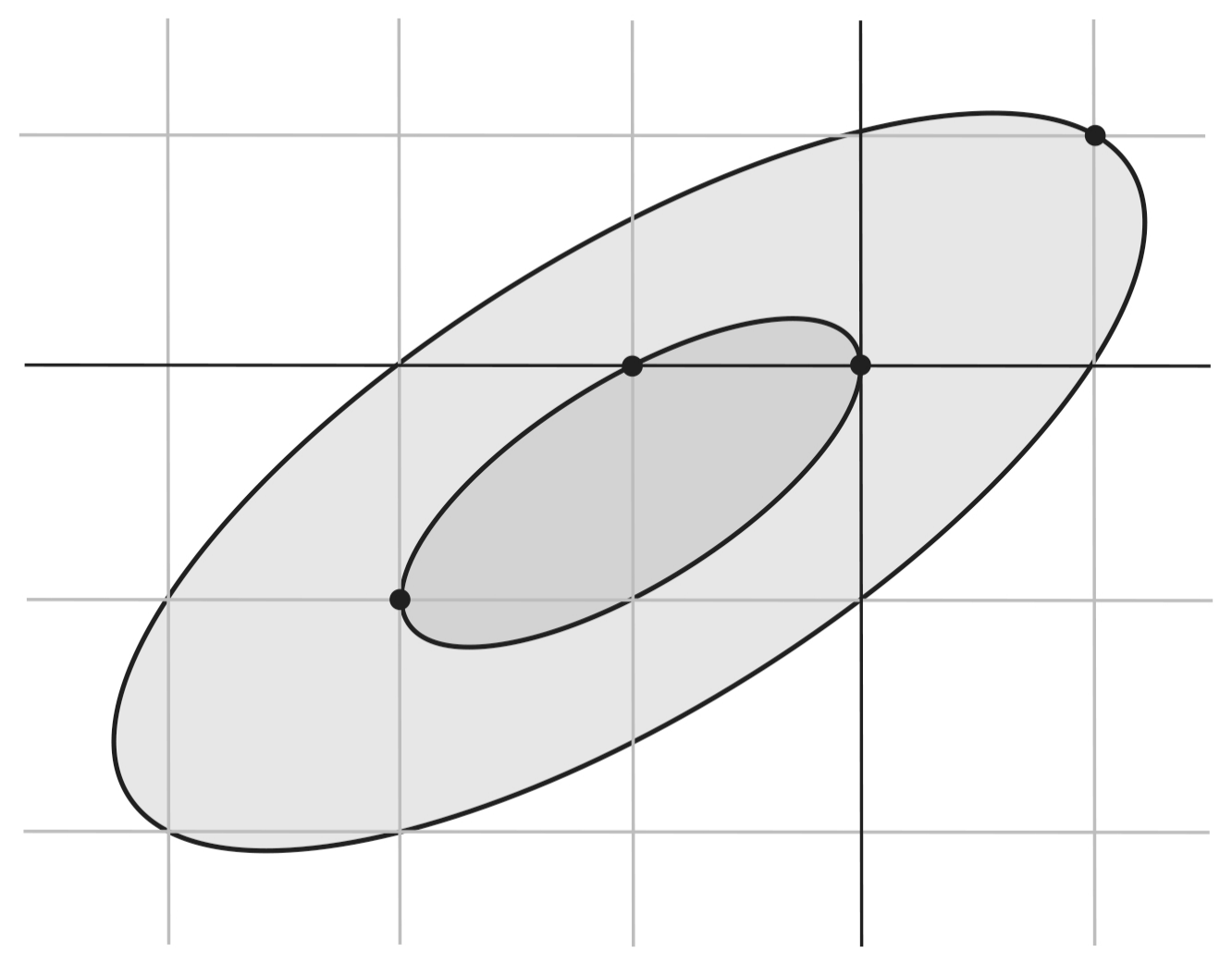}
{\Small
\put(-95,58){$S_0$}
\put(-130,31){$S_1$}
\put(-118,51){$1$}
\put(-93,82){$-1$}
\put(-48,82){$-1$}
\put(-16,113){$1$}
}}
\hspace{3cm}
{\begin{tikzpicture}
\node (0) at (0,0) {$\bullet$};
\node (1) at (0,1) {$\bullet$};
\node (2) at (0,2) {$\bullet$};
\node (3) at (0,3) {$\bullet$};
\node at (0,3.75) {$\vdots$};
\node[right] at (2.05,3.75) {$\vdots$};
\draw (0,0) -- (0,3);

\draw[dashed] (0,0) --++ (2,0);
\draw[dashed] (0,1) --++ (2,0);
\draw[dashed] (0,2) --++ (2,0);
\draw[dashed] (0,3) --++ (2,0);

\node[right] at (2,0) {$0$};
\node[right] at (2,1) {$2$};
\node[right] at (2,2) {$4$};
\node[right] at (2,3) {$6$};

\node[left] at (-.2,0) {$-1$};
\node[left] at (-.2,1) {$0$};
\node[left] at (-.2,2) {$0$};
\node[left] at (-.2,3) {$0$};

\end{tikzpicture}
}
    \caption{Left: the sublevel sets $S_0$  and $S_1$ from Example \ref{ex:sublevel set example}. The four contributing lattice points and their $F_{\Gamma,k}$-weights are indicated. Right: the corresponding graded root, with the weights from Equation \eqref{eq:sum over connected component} written to the left of each vertex. Recall from grading conventions \ref{def:normalizating graded root} that the numbers $0,2,4,6$ indicate the Heegaard Floer gradings. }
    \label{fig:sublevel set and weights example}
\end{figure}

\begin{exmp}
\label{ex:sublevel set example}
We end this section with an example of the above weights $F_{\Gamma,k}$. See Figure \ref{fig:sublevel set and weights example} for a summary. For the purpose of a $2$-dimensional illustration, consider the plumbing representation $\Gamma$ for $S^3$ with two vertices of weight $m_1=-1$ and $m_2=-2$, shown below.
\begin{equation*}
    \begin{tikzpicture}
    \node at (0,0) {$\bullet$};
    \node at (1.5,0) {$\bullet$};
    \node[above] at (0,0) {$-1$};
    \node[above] at (1.5,0) {$-2$};
    \draw (0,0) --++ (1.5,0) ;
    \end{tikzpicture}
\end{equation*}We pick the $\spinc$ representative $k=(-1,0)$, so that for $(x,y)\in \Z^2$ we have 
\begin{equation*}
    \chi^{}_k(x,y) = (x^2+2y^2-2xy+x)/2, \hskip1em
    F_{\Gamma,k}(x,y) = F_1(-2x+2y-1)F_1(2x-4y+1). 
\end{equation*}
By the formula for $F_1$ in \eqref{eq:F_1 and F_0}, we see that $F_{\Gamma,k}(1,1) = F_{\Gamma,k}(-2,-1)  = 1$,  $F_{\Gamma,k}(-1,0) = F_{\Gamma,k}(0,0) = -1$, and all other lattice points have weight zero and thus do not contribute.

One may verify that the minimum value of $\chi^{}_k$ is zero and that all sublevel sets are connected. The sublevel set $S_0$ contains $(-2,-1), (-1,0)$, and $(0,0)$, but not $(1,1)$, and $S_1$ contains $(1,1)$. The weight assigned to the vertex in the graded root corresponding to $S_0$ and $S_j$ for $j\geq 1$ is then $-1$ and $0$, respectively. Compare with Figure \ref{fig:wgroot -1}, which specializes to the above discussion at $q=t=1$. 
\end{exmp}

\section{The invariant} \label{sec: The invariant}

This section introduces the main construction of this paper, a refinement of the weights from equation \eqref{eq:sum over connected component} in the form of a collection of two-variable Laurent polynomials. 
Section \ref{sec: invariance} shows that the resulting weighted graded root is a $3$-manifold invariant.

\subsection{Refined weights} 
We start by establishing the following notation. 

\begin{notn} \label{notn}
For $k,x \in \Z^s$, define
\begin{align*}
    \normalization_k &= -\frac{(k-Mu)^2 + 3s + \sum_v m_v}{4},& \tnormalization_k &=\frac{k\cdot u-\langle u,u\rangle}{2},\\ 
    \eps_k(x) &= \normalization_k + 2\chi^{}_k(x) + \langle x,u\rangle,&
    \theta_{k}(x) &= \tnormalization_k+\langle x,u\rangle,
\end{align*}
where the notation $(-)^2$ is the same as in Remark \ref{rem: k^2 meaning}. Also, recall that $Mu = m+\delta$, so $\langle x,u\rangle = x\cdot (m+ \delta)$ and $\langle u, u\rangle = \sum\limits_{v}(m_v+\delta_v)$.

The terms $\normalization_k$ and $\tnormalization_k$ are overall normalizations\footnote{We thank the referee and  Sunghyuk Park for pointing out the above normalization $\tnormalization_k$.} used to eliminate dependence on the choice of $\spinc$ representative. In general, $\normalization_k$ is a rational number and is similar in form to the $d$-invariant from Heegaard Floer homology (see Remark \ref{rem: d inv}). On the other hand, $\tnormalization_k$ is always an integer because $k$ is characteristic. Also, note $\tnormalization_k = \chi^{}_{-k}(u)$.
\end{notn}

\begin{defn}
\label{def:weighted graded root k}
Let $k\in m + 2\Z^s$ be a $\spinc$ representative and let $F=\{F_n\}_{n\geq 0}$ be an admissible family valued in a commutative ring $\Ring$. To each vertex of the graded root $\groot{k}$ we assign a weight valued in $q^{\normalization_k}\cdot\Ring[q^{\pm 1}, t^{\pm 1}]$ as follows. For a vertex represented by a connected component $C$ in some sublevel set, let $L(C) = C \cap \Z^s$ denote its lattice points. Set 
\begin{equation}
    \label{eq:laurent poly enhancement}
    P_{F,k}(C) = \sum\limits_{x\in L(C)}F_{\Gamma, k}(x)q^{\eps_k(x)}t^{\theta_k(x)}, 
\end{equation}
and let $\wgrootF{k}$ denote the graded root $\groot{k}$ with these weights. We will often omit the reference to  $F$ by writing $P_k$ instead of $P_{F,k}$. Note that specializing $q=t=1$ recovers the weights in equation \eqref{eq:sum over connected component}.
\end{defn}

The above weights can be interpreted geometrically as follows. For $n\in \Z$, the coefficient of $t^{n+\tnormalization_k}$ in $P_{k}(C)$ is given by summing $F_{\Gamma,k}(x)q^{\normalization_k + 2\chi^{}_k(x) + n}$ over all $x\in \Z^s$ which lie on the intersection of $C$ with the hyperplane $\{y\in \R^s \mid \langle y,u\rangle = n \}$.  

Let us verify that the weights $P_{k}(C)$ are compatible with the isomorphisms (\ref{eq:GR for different reps}) relating graded roots for different representatives of $[k]$. 

\begin{lem}
\label{lem:equations for different spinc reps}
Let $k, k'=k+2My \in m+2\Z^s$ be two representatives for $[k]\in \spinc (Y)$. Then $\eps_{k'}(x) = \eps_k(x+y)$ and $\theta_{k'}(x) = \theta_k(x+y)$ for all $x\in \Z^s$.
\end{lem}

\begin{proof}
First, we show that $\eps_{k'}(x) = \eps_k(x+y)$. Note that
\[
(k'-Mu)^2 = (k-Mu)^2 +4 y\cdot k - 4 \langle y,u\rangle + 4 \langle y,y\rangle, 
\]
which implies 
\begin{align*}
    \normalization_{k'} &= \normalization_k - y\cdot k +\langle y,u\rangle - \langle y,y\rangle  \\
    &= \normalization_k + 2\chi^{}_{k}(y) + \langle y,u\rangle. 
\end{align*}
The equality $\eps_{k'}(x) = \eps_k(x+y)$ now follows from the above equality and equation \eqref{eq:chi for different reps}.
Next,
\begin{align*}
    \tnormalization_{k'} = \frac{(k+2My)\cdot u-\langle u,u\rangle}{2} = \tnormalization_k +\langle y,u\rangle
\end{align*}
Hence, $\theta_{k'}(x) = \tnormalization_k+\langle x+y, u\rangle = \theta_k(x+y)$.
\end{proof}

\begin{prop}
\label{prop:different reps for weighted graded root}
Let $k\in m + 2\Z^s$ and let $k' = k + 2My$ for some $y\in \Z^s$. Then 
\[
\wgroot{k'} \cong \wgroot{k}[-\chi^{}_{k}(y)]. 
\]
\end{prop}

\begin{proof}
Recall that the isomorphism $\groot{k} \cong \groot{k'}[-\chi^{}_{k}(y)]$ of graded roots in equation \eqref{eq:GR for different reps} is induced by the translation $T(x) = x+y$. A lattice point $x\in L(C)$ contributes the summand
\[
F_{\Gamma, k'}(x)q^{\eps_{k'}(x)}t^{\theta_{k'}(x)}
\]
to $P_{F,k'}(C)$, while $T(x)$ contributes 
\[
F_{\Gamma, k}(x+y)q^{\eps_k(x+y)}t^{\theta_k(x+y)}
\]
to $P_{F,k}(T(C))$. Lemma \ref{lem:equations for different spinc reps} and equation \eqref{eq:F different representatives} imply that
\[
F_{\Gamma, k'}(x)q^{\eps_{k'}(x)}t^{\theta_{k'}(x)} = F_{\Gamma, k}(x+y)q^{\eps_k(x+y)}t^{\theta_k(x+y)},
\]
which completes the proof. 
\end{proof}

\begin{defn}
\label{def:weighted graded root [k]}
Set $\wgroot{[k]}$ to be $\groot{[k]}$, as in Definition \ref{def:normalizating graded root},  equipped with the weights $P_k$ for some $k\in [k]$. 
\end{defn}

\noindent Lemma \ref{lem:equations for different spinc reps} and Proposition \ref{prop:different reps for weighted graded root} guarantee that $\wgroot{[k]}$ does not depend on $k\in [k]$.

\subsection{Invariance} \label{sec: invariance}

In this section we prove invariance of $\wgroot{[k]}$ under the two Neumann moves shown in Figure \ref{fig:Neumann moves}. This establishes Theorem \ref{main thm}, which is restated as Theorem \ref{thm:invariance} below using a more detailed notation. In what follows, $\Gamma$ is a negative definite plumbing tree with $s$ vertices, and $\Gamma'$ is a plumbing tree with $s'=s+1$ vertices obtained from $\Gamma$ by one of the type (a) or (b) moves. We will use the conventions established in Notation \ref{notation type a and b}, as well as the following additional notation for the two moves.

{\bf Type (a):} The intersection form for $\Gamma'$ is given by $M' = \til{M} + A$ where 
\begin{align}
\label{eq:M' for type a} 
    \til{M} =
\begin{pmatrix}
0 & 0 \\
0 & M
\end{pmatrix}
\hskip2em 
A = 
{\scriptsize
\begin{pmatrix}
-1 & \phantom{-}1 & \phantom{-}1 & 0 & \cdots & 0 \\
\phantom{-}1  & -1  & -1  &  0 &  \cdots      &  0 \\
\phantom{-}1 &   -1   &   -1 & 0 & \cdots & 0 \\
\phantom{-}0 & \phantom{-}0 & \phantom{-}0 & 0 & \cdots & 0 \\
\phantom{-}\vdots &  & & & & \vdots \\
\phantom{-}0 & \phantom{-}\cdots & \phantom{-}\cdots & \cdots & \cdots & 0
\end{pmatrix}
}
\end{align}

As in \cite[Proposition 3.4.2]{Nem_Lattice_cohomology}, define the projection $\pi_{*}:\Z^{s+1}\to \Z^s$ by 
\begin{align}
\label{eq:pi_* map} 
    \pi_{*}(x_{0}, x_{1}, \ldots, x_{s}) = (x_{1}, \ldots, x_{s})
\end{align}
and the inclusion $\pi^*:\Z^s\to\Z^{s+1}$ by
\begin{align}
\label{eq:pi^* map type a}
    \pi^*(x_{1}, x_2, \ldots, x_{s}) = (x_{1}+x_{2}, x_{1}, x_{2}, \ldots, x_{s}).
\end{align}

{\bf Type (b):} The intersection form for $\Gamma'$ is given by $M' = \til{M} + A$ where
\begin{align} \label{M prime type b}
    \til{M}=\begin{pmatrix}
0 & 0 \\
0 & M
\end{pmatrix}
\hskip2em 
A = 
{\scriptsize
\begin{pmatrix}
-1 & \phantom{-}1 & 0 & \cdots & 0 \\
\phantom{-}1  & -1  &  0 &  \cdots      &  0 \\
\phantom{-}0 & \phantom{-}0 & 0 & \cdots & 0 \\
\phantom{-}\vdots &   & & & \vdots \\
\phantom{-}0 & \phantom{-}\cdots & \cdots & \cdots & 0
\end{pmatrix}
}
\end{align}

As in the type (a) case, let $\pi_{*}:\Z^{s+1}\to\Z^s$ denote the projection $\pi_*(x_0,x_{1}, \ldots, x_{s}) = (x_{1}, x_{2}, \ldots, x_{s})$, and define two inclusions $\pi^*, \rho^* : \Z^s \to \Z^{s+1}$ by
\begin{align}
    \pi^*(x_{1}, \ldots, x_{s}) &= (x_{1}, x_{1}, \ldots, x_{s}), \label{eq:pi^* map type b} \\
    \rho^* (x_{1}, \ldots, x_{s}) &= (x_{1}-1, x_{1},  \ldots, x_{s}). \label{eq:p map type b}
\end{align}

With the notation in place, we now record several results regarding various contributions to the Laurent polynomial weights. 

\begin{lem} \label{lem: exponent of normalization term}
Let $\Gamma$, $\Gamma'$ be negative definite plumbings related by a Neumann move as above. Let $k\in m+2\Z^s$ be a $\spinc$ representative, and let $k'\in m'+2\Z^{s+1} $ denote the associated representative, as in equations \eqref{eq:k' for type a} and \eqref{eq:k' for type b}. Then $\normalization_k = \normalization_{k'}$ and $\tnormalization_k = \tnormalization_{k'}$. 
\end{lem}
\begin{proof}
Let $u'=(1,u)$.
\\
\underline{Type (a):} First we show $\normalization_k = \normalization_{k'}$. Observe that $3s + \sum_v m_v = 3s' + \sum_v m_v'$, so it remains to verify that 
\[
(k'-M'u')^2 = (k-Mu)^2. 
\]
Expressions for $M'$ and $k'$ are given in equations \eqref{eq:M' for type a} and \eqref{eq:k' for type a}. Note that 
\[
k' - M'u' =  (0,k-Mu).
\]
Let $y=(y_1, y_2, \ldots, y_{s})\in \Q^{s}$  be such that $k-Mu = My$. Then $M'(y_1 + y_2, y) = (0,k-Mu)$, and it follows that $(k'-M'u')^2 = (k-Mu)^2$. Next, we show $\tnormalization_k = \tnormalization_{k'}$.
\begin{align*}
    \tnormalization_{k'} = \frac{[(0,k)+(1,-1,-1,0,\ldots, 0)]\cdot (1,u)-\langle (1,u), (1,u)\rangle}{2}&=\frac{k\cdot u -1-(1,u)^t(\widetilde{M}+A)(1,u)}{2}\\
    &=\frac{k\cdot u -1 -u^tMu +1}{2} = \tnormalization_k
\end{align*}
\underline{Type (b):} First we show $\normalization_k = \normalization_{k'}$. $M'$ and $k'$ are given in equations \eqref{M prime type b}, \eqref{eq:k' for type b}. Note 
$$k'-M'u'= (-1,1, 0,\ldots,0)+(0, k-Mu).
$$

Denote these two summands by $w:=(-1,1, 0,\ldots,0)$ and $\widetilde k:= (0, k-Mu)$.
We claim that
\begin{equation} \label{eq: claim}
(\widetilde k)^t(M')^{-1}\widetilde k = (k-Mu)^t M^{-1}(k-Mu). 
\end{equation}

To prove the claim, let $y=M^{-1}(k-Mu)$ or equivalently $k-Mu=My$. Then one checks that $\widetilde k=M'(y_1, y)$, where $y_1$ is the first coordinate of $y$. Thus the left-hand side of (\ref{eq: claim}) equals $(0, k-Mu)^t(y_1, y)$ which equals the right-hand side of (\ref{eq: claim}), verifying the claim.
Expanding linearly, consider

$$(k'-M'u')^t(M')^{-1}(k'-M'u') - (k-Mu)^t M^{-1}(k-Mu)=
$$

\begin{equation} \label{eq: three terms}
w^t(M')^{-1}w + 2(\widetilde k)^t(M')^{-1} w.
\end{equation}

It follows from equations (\ref{M prime type b}) that $(M')^{-1} w = (1,0,\ldots,0),$ thus the first term in (\ref{eq: three terms}) equals $-1$ and the second term is zero. Under the (b) move, we have $3s+\sum_{v}m_{v} = 3s' + \sum_v m_v' -1$, precisely offsetting the change in $(k-Mu)^2$ computed above. Therefore, $\normalization_k = \normalization_{k'}$. Next, we show $\tnormalization_k = \tnormalization_{k'}$.
\begin{align*}
    \tnormalization_{k'} = \frac{[(0,k)+(-1,1,0,\ldots, 0)]\cdot (1,u)-\langle (1,u), (1,u)\rangle}{2}&=\frac{k\cdot u -(1,u)^t(\widetilde{M}+A)(1,u)}{2}\\
    &=\frac{k\cdot u -u^tMu }{2} = \tnormalization_k
\end{align*}
\end{proof}

For the type (a) move in the following lemma, recall the function $\pi^*$ from equation \eqref{eq:pi^* map type a}. 
For the type (b) move, recall the functions $\pi^*$, and $\rho^*$ from equations \eqref{eq:pi^* map type b} and \eqref{eq:p map type b}. 

\begin{lem}
\label{lem:pi* maps and u}
Let $\Gamma, \Gamma', k$, and $k'$ be as in the statement of Lemma \ref{lem: exponent of normalization term}. 
\begin{enumerate}
    \item In the type (a) case, for any $x\in \Z^{s}$,    
\[
        \langle x , u \rangle  =\langle  \pi^*(x) , u'\rangle.
\]
    \item  In the type (b) case, for any $x\in \Z^s$,
\[
       \langle x ,u\rangle  =\langle \pi^*(x) , u'\rangle = \langle \rho^*(x),u'\rangle.
\]
\end{enumerate}
\end{lem}

\begin{proof}
Recall that $\langle x,u\rangle = x\cdot (m+\delta)$. For item (1), $m' + \delta' = (1,-1,-1, 0, \ldots, 0) + (0, m+\delta)$. Then 
\begin{align*}
    \pi^*(x)\cdot (m' +\delta') &= (x_{1}+x_{2}, x_{1}, x_{2}, \ldots, x_{s})\cdot [(1,-1,-1,0, \ldots, 0) + (0, m+\delta)]\\
    &= x\cdot (m+\delta).
\end{align*}
For item (2), $m' + \delta' = (0,m + \delta)$, and the desired equality follows. 
\end{proof}

\begin{lem}
\label{lem:F functions and pi maps}
Let $\Gamma, \Gamma', k$, and $k'$ be as in the statement of Lemma \ref{lem: exponent of normalization term}. 
\begin{enumerate}
    \item For the type (a) move, $F_{\Gamma, k}(x) =  F_{\Gamma', k'}(\pi^*(x))
    $ for all $x\in \Z^s$, 
    \item For the type (b) move,  $F_{\Gamma, k}(x) =  F_{\Gamma',k'}(\pi^*(x)) + F_{\Gamma',k'}(\rho^*(x))
    $ for all $x\in \Z^s$.

\end{enumerate}
\end{lem}

\begin{proof}
To show item (1), first note that 
\[
    2M'\pi^*(x) + k' - M'u' =(0, 2Mx+ k -Mu), 
\]
so, using property \ref{item:A1}, $F_{\Gamma', k'}(\pi^*(x)) = F_{2}(0)\prod\limits_{i=1}^{s} F_{\delta_{i}}((2Mx + k -Mu)_{v_{i}})=F_{\Gamma, k}(x)$.

We now verify item (2). Observe that
\begin{align*}
2M'\pi^*(x) + k' - M'u' &= (0,2Mx+k-Mu) + (-1,1,0,\ldots, 0), \\
2M'\rho^*(x) + k' - M'u' &= (0, 2Mx +k-Mu)+ (1,-1, 0, \ldots, 0). 
\end{align*}
Introduce the notation
\begin{align*}
    r = (2Mx + k - Mu)_{v_1},\hspace{2em} \overline{r} = \prod_{i=2}^s F_{\delta_i}(2Mx + k - Mu)_{v_i},
\end{align*}
and recall from equation \eqref{eq:F_1 and F_0} that $F_1(\pm 1) = \mp 1$.  Then we have 
\begin{align*}
    F_{\Gamma',k'}(\pi^*(x)) &= F_{\delta_1 +1}(r + 1) \cdot \overline{r}, \\
    F_{\Gamma',k'}(\rho^*(x)) &= -F_{\delta_1 +1 }(r-1) \cdot \overline{r}, \\
    F_{\Gamma,k}(x) &= F_{\delta_1}(r) \cdot \overline{r},
\end{align*}
and the desired equality follows from property \ref{item:A2}. 
\end{proof}

We are in a position to prove our main result:

\begin{thm}
\label{thm:invariance} 
For any admissible family of functions $F$, the weighted graded root $\wgroot{[k]}$ is an invariant of the $3$-manifold $Y(\Gamma)$ equipped with the $\spinc$ structure $[k]$. 
\end{thm}

\begin{proof}
We will demonstrate an isomorphism $
\wgroot{k} \cong \wgroot{k'}$ of weighted graded roots when $\Gamma, \Gamma', k$, and $k'$ are as in the statement of Lemma \ref{lem: exponent of normalization term}. 

For each of the two moves we first give an explicit isomorphism of graded roots $\groot{k}\cong \groot{k'}$, following the proofs of \cite[Proposition 4.6]{Nem_On_the} and \cite[Proposition 3.4.2]{Nem_Lattice_cohomology}. We then show that this isomorphism respects our Laurent polynomial weights.

We begin with the type (a) move. Recall the functions $\pi_*$ and $\pi^*$ from \eqref{eq:pi_* map} and \eqref{eq:pi^* map type a}, and that $M' = \widetilde{M} + A$, as in equation \eqref{eq:M' for type a}.

For $x'= (x_0, x_1, x_2, \ldots, x_s) \in \Z^{s+1}$, we have
\[
(x')^t A x' = - (x_0 - x_1 - x_2)^2,
\]
so $(x')^t M' x' = \pi_*(x')^t M \pi_*(x') - (x_0 - x_1 - x_2)^2$. It is then straightforward to verify that 
\[
 \chi^{}_{k'}(x')  = \chi^{}_{k}(\pi_{*}(x')) + \frac{1}{2}[(x_{0}-x_{1}-x_{2})(x_{0}-x_{1}-x_{2}-1)].
\]
In particular, substituting $x'=\pi^*(x)$ for $x\in \Z^s$, this implies 
\begin{equation}
\label{eq:pi^* commutes with chi type a}
\chi^{}_{k'} \circ \pi^* = \chi^{}_{k},
\end{equation}
so $\pi^*$ induces an inclusion $  \chi_{k}^{-1}((-\infty, j])\hookrightarrow \chi_{k'}^{-1}((-\infty, j])$ of sublevel sets.
As in the proof of \cite[Proposition 3.4.2]{Nem_Lattice_cohomology}, $\pi^*$ also induces a bijection, denoted $\tilde{\pi}^*$, on connected components in these sublevel sets. The isomorphism $\groot{k} \cong \groot{k'}$ of graded roots is given by $\tilde{\pi}^*$, sending a connected component $C$ to the connected component $\til{\pi}^*(C)$ that contains $\pi^*(C)$.  

Fix a connected component $C$ in some sublevel set of $\chi^{}_k$. We will now show that
\[
P_k(C) = P_{k'}(\tilde{\pi}^*(C)).
\]
The term on the right-hand side above is a sum over contributions from lattice points in the component $\tilde{\pi}^*(C)$, which contains all the lattice points in $\pi^*(C)$, but is in general strictly bigger. As we shall now see, only lattice points in $\pi^*(C)$ contribute. We have
\[
(2M'x'+ k'-M'u')_{v_0'} = -2(x_0 - x_1 - x_2).
\]
Since $\delta_{v_0'}= 2$, property \ref{item:A1} implies that $F_{\Gamma',k'}(x') = 0$ unless $x_0 = x_1 + x_2$, so 
\[
P_{k'}(\tilde{\pi}^*(C)) =  \sum_{x'\in \pi^*(L(C))} F_{\Gamma', k'}(x') q^{\eps_{k'}(x')} t^{\theta_{k'}(x')}.
\]
Therefore, it suffices to show 
\[
F_{\Gamma, k}(x) q^{\eps_k(x)} t^{\theta_k(x)} 
=
F_{\Gamma', k'}(\pi^*(x)) q^{\eps_{k'}(\pi^*(x))} t^{\theta_{k'}(\pi^*(x))}.
\]
for all $x\in C$. Equation \eqref{eq:pi^* commutes with chi type a}, Lemma \ref{lem: exponent of normalization term}, and item (1) of Lemma \ref{lem:pi* maps and u} guarantee that the powers of $q$ and $t$ are equal, and $F_{\Gamma, k}(x) = F_{\Gamma', k'}(\pi^*(x))$ by Lemma \ref{lem:F functions and pi maps} (1). This concludes the proof of the type (a) move.

We now address the type (b) move. Recall the functions $\pi_*, \pi^*$, and $\rho^*$ from equations \eqref{eq:pi_* map}, \eqref{eq:pi^* map type b}, \eqref{eq:p map type b}, and that $M'= \widetilde{M} + A$ as in equation \eqref{M prime type b}. For $x'\in \Z^{s+1}$, we have
\[
(x')^t A x' = - (x_0 - x_1)^2,
\]
so $(x')^t M' x' = \pi_*(x')^t M \pi_*(x') - (x_0 - x_1)^2$. It is then easy to see that 
\[
\chi^{}_{k'}(x') = \chi^{}_{k}(\pi_*(x')) +\frac{1}{2}(x_{0}-x_{1})(x_{0}-x_{1}+ 1),
\]
which implies 
\begin{equation}
      \label{eq:pi^* commutes with chi type b}
    \chi^{}_{k'} \circ \pi^*  = \chi^{}_k = \chi^{}_{k'} \circ \rho^*.
\end{equation}

Thus both $\pi^*$ and $\rho^*$ induce inclusions $  \chi_{k}^{-1}((-\infty, j])\hookrightarrow \chi_{k'}^{-1}((-\infty, j])$ of sublevel sets. As in the type (a) case above, $\pi^*$ also induces a bijection $\tilde{\pi}^*$ between connected components of each sublevel set, and the isomorphism of graded roots $\groot{k} \cong \groot{k'}$ is given by $\tilde{\pi}^*$. 

To complete the proof we check that 
\[
P_k(C) = P_{k'}(\tilde{\pi}^*(C))
\]
for every connected component $C$ of each sublevel set of $\chi^{}_k$. As in the type (a) case, we will now see that only a particular subset of lattice points in $\tilde{\pi}^*(C)$ contribute to $P_{k'}(\tilde{\pi}^*(C))$. To begin, note
\[
(2M'x' + k' - M'u')_{v_0'} = -2(x_0 - x_1) - 1.
\]
Since $\delta_{v_0'} =1$, the formula for $F_1$ from equation \eqref{eq:F_1 and F_0} implies that $F_{\Gamma',k'}(x') = 0$ unless $-2(x_0 - x_1) - 1 = \pm 1$, or, equivalently, unless $x' = \pi^*(x)$ or $x' = \rho^*(x)$ for some $x\in \Z^s$. Observe that $\pi^*(x) - \rho^* (x) = (1,0, \ldots, 0)$, so $\pi^*(x)$ and $\rho^*(x)$ are in the same component of $\chi_{k'}^{-1}(j)$. It follows that
\[
P_{k'}(\tilde{\pi}^*(C)) = 
\sum_{x'\in \pi^*(L(C))} F_{\Gamma', k'}(x') q^{\eps_{k'}(x')} t^{\theta_{k'}(x')}
+
\sum_{x'\in \rho^*(L(C))} F_{\Gamma', k'}(x') q^{\eps_{k'}(x')} t^{\theta_{k'}(x')}.
\]
To complete the proof, we have
\[
F_{\Gamma, k}(x) q^{\eps_k(x)} t^{\theta_k(x)} 
=
  F_{\Gamma', k'}(\pi^*(x)) q^{\eps_{k'}(\pi^*(x))} t^{\theta_{k'}(\pi^*(x))} 
 +   F_{\Gamma', k'}(\rho^*(x)) q^{\eps_{k'}(\rho^*(x))} t^{\theta_{k'}(\rho^*(x))}
\]
by combining Lemma \ref{lem: exponent of normalization term}, equation \eqref{eq:pi^* commutes with chi type b}, item (2) of Lemma \ref{lem:pi* maps and u}, and item (2) of Lemma \ref{lem:F functions and pi maps}. 
\end{proof}

\section{The two-variable series}
\label{sec:two variable series}

In this section we extract a two-variable series from $\wgroot{[k]}$ by taking a limit (in an appropriate sense) of the weights $P_k(C)$. Theorem \ref{thm:convergence} shows that this limiting procedure yields a well-defined invariant of $(Y(\Gamma), [k])$.  Throughout this section  some $\Ring$-valued admissible family $F$ will be fixed, and references to it will be omitted for brevity of notation.

We first establish some preliminary notions. For a commutative ring $\Ring$, denote by $\Ring[q^{-1},q]]$ the ring of Laurent series in $q$ and by $\Ringt$ the set of Laurent series in $q$ whose coefficients are Laurent polynomials in $t$, 
\[
\Ringt = \left(\Ring[t^{\pm 1}] \right)[q^{-1},q]].
\]
Given $\normalization\in \Q$, $\tnormalization\in\Z$, $f\in  q^\normalization \cdot \Ringt $, and $i,j \in \Z$, let $[f]_{i,j} \in \Ring$
be the coefficient of $q^{\normalization+i} t^{\tnormalization+j}$ in $f$.

\begin{defn}
\label{def:stabilization}
We say a sequence $f_1, f_2, \ldots  \in q^\normalization \cdot \Ring [q^{\pm 1}, t^{\pm 1}]$ \emph{stabilizes} if for all $i,j\in \Z$, the sequence of coefficients $([f_1]_{i,j}, ([f_2]_{i,j}, \ldots )$ is eventually constant. For such a sequence, the \emph{limit} $f$ is the bi-infinite series in $q,t$ defined
by setting $[f]_{i,j}$ to be the limit of $[f_n]_{i,j}$ as $n\to \infty$. 
\end{defn}

As stated in the definition, the limit of a stabilizing sequence in general is a bi-infinite series in $q,t$. In Theorem \ref{thm:convergence} below, the limit is claimed to be an element of $q^\normalization \cdot \Ringt$. In addition to proving that the sequence $(f_n)$ stabilizes, this will be shown by establishing that \begin{enumerate}[label=(\roman*)]
\item
there exists $i_0\in \Z$ such that $[f_n]_{i,j} = 0$ for all $n\geq 0$, $j\in \Z$, and $i\leq i_0$, and 
\item for any fixed $i$, the set of $j$ such that $[f_n]_{i,j} \neq 0$ is bounded.
\end{enumerate}

Returning to weighted graded roots, fix a negative definite plumbing tree $\Gamma$ and $\spinc$ representative $k\in m+2\Z^s$. Consider the weighted graded root $\wgroot{k}$, as given in Definition \ref{def:weighted graded root k}. For $n\in \Z$, let 
\begin{equation}
\label{eq:P_k^n}
P_k^n := \sum_{C\in \chi^{-1}_{k}(n)} P_k(C) \in q^{\normalization_k} \cdot \Ring[q^{\pm 1}, t^{\pm 1}] 
\end{equation}
denote the sum of the Laurent polynomial weights over vertices $C$ of $R_k$ in $\chi^{}_{k}$-grading $n$.
Recall that $\chi^{}_k$ is bounded below by some $n_0\in \Z$, and consider the sequence $(P_k^{n_0}, P_k^{n_0 + 1}, P_k^{n_0 + 2}, \ldots )$. 

\begin{rem}
Note that $P_k^n$ is the sum of $F_{\Gamma,k}(x)q^{\eps_k(x)}t^{\theta_k(x)}$ over all lattice points $x$ in the entire $n$-sublevel set of $\chi^{}_{k}$. Moreover, since there is only one connected component for large enough $n$, one may just as well start the sequence at a sufficiently high $\chi^{}_{k}$-grading, making the sum in \eqref{eq:P_k^n} be given by a single $P_k(C)$. 
\end{rem}

\begin{thm}
\label{thm:convergence}
The sequence $\left(P_k^{n_0}, P_k^{n_0 + 1}, P_k^{n_0 + 2}, \ldots\right)$ stabilizes to an element of $q^{\normalization_k} \cdot  \Ringt$. Its limit, which we denote $P_{[k]}^\infty$, is an invariant of the $3$-manifold $Y(\Gamma)$ equipped with the $\spinc$ structure $[k]$.
\end{thm}

\begin{proof}
For $n\in \Z$, define
\begin{align*}
    \mathcal{S}_n &= \{ x\in \Z^s \mid \chi^{}_k(x) \leq n \}, &
   \til{\mathcal{S}}_n &= \{ x\in \Z^s \mid 2\chi^{}_k(x) + \langle x, u\rangle  \leq n \}, \\
   \partial \til{\mathcal{S}}_n &= \{ x\in \Z^s \mid 2\chi^{}_k(x) + \langle x, u\rangle =n \},
     & \mathcal{A}_n &= \{ x\in \Z^s \mid \langle x,u\rangle = n\}.
\end{align*} 
By definition, 
\[
P_k^n = \sum_{x\in \mathcal{S}_n} F_{\Gamma,k}(x)q^{\eps_k(x)} t^{{\theta_k(x)}}.
\]
It follows from Notation \ref{notn} that for fixed $i,j\in \Z$, the coefficient of $q^{\normalization_k + i} t^{\tnormalization_k+j}$ in $P_k^n$ is equal to 
\begin{equation}
\label{eq:coefficient fixed n i and j}
\sum_{x} F_{\Gamma,k}(x),
\end{equation}
where the sum is over $x\in \mathcal{S}_n \cap \partial\til{\mathcal{S}}_i \cap \mathcal{A}_j$. Both $(\mathcal{S}_n)$, $(\til{\mathcal{S}}_n)$ are sequences of nested finite sets whose union is $\Z^s$. Hence for a fixed $i$ there exists $N$ such that $\sublevel_i \subset \sublevel_n$ for all $n\geq N$. Then for $n\geq N$, we have 
\[
\mathcal{S}_n \cap \partial\til{\mathcal{S}}_i \cap \mathcal{A}_j = \partial\til{\mathcal{S}}_i \cap \mathcal{A}_j,
\]
so that the sum in equation \eqref{eq:coefficient fixed n i and j} is independent of $n$ for $n$ sufficiently large. See Figure \ref{fig:sublevel sets} for an illustration when $s=2$. This verifies stabilization of the sequence. 

As discussed after Definition \ref{def:stabilization}, we will check two conditions (i), (ii) ensuring that the limit is an element of $q^\normalization \cdot \Ringt$.
The condition (i) follows from the fact that the $q$-powers in the $P_k^{n}$ are given by $\eps_k$, which is bounded below. To check (ii), observe that for a fixed $i$, the exponent of $t$ is given by $\tnormalization_k+\langle x, u\rangle$ which is bounded on the set $\til{\mathcal{S}}_i$.
\begin{figure}
    \centering
    \begin{tikzpicture}
\draw[thick, rotate=-10 , fill=gray!35] (0,0) ellipse (1 and .75);
\draw[rotate = 20, thick] (-.5,.6) ellipse (.5 and .75);
\draw[thick] (-2.5,1.25) -- (2.5,-1.25);
\node at (0-.15,0+.075) {$\bullet$};

\begin{scope}[shift= {(7,0)}]
\draw[thick, rotate=-10 , fill=gray!35] (0,0) ellipse (2 and 1.3);
\draw[rotate = 20, thick] (-.5,.6) ellipse (.5 and .75);
\draw[thick] (-2.5,1.25) -- (2.5,-1.25);
\node at (0-.15,0+.075)  {$\bullet$};
\node at (0-1.2,0+.6) {$\bullet$};
\node at (0,-1.25) {};
\end{scope}
\end{tikzpicture}
    \caption{A schematic depiction of stabilization. Left: $\sublevel_n$ is represented by the shaded ellipse, $\partial\til{\sublevel}_i$ is the unshaded ellipse, and the hyperplane $\mathcal{A}_j$ is the straight line (the actual sets are discrete subsets of the illustration). The intersection $\mathcal{S}_n \cap \partial\til{\mathcal{S}}_i \cap \mathcal{A}_j$ is marked by a dot. Right: taking $n>>0$ ensures $\partial\til{\sublevel}_i \subset \sublevel_n$, and the intersection (two dots) is the same for all sufficiently large $n$.} 
    \label{fig:sublevel sets}
\end{figure}
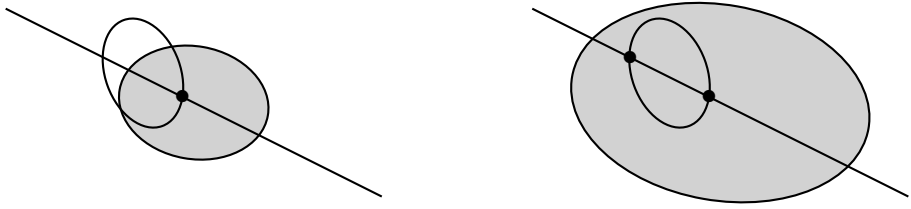

That $P_{[k]}^\infty$ is an invariant of $(Y(\Gamma),[k])$ follows from Proposition \ref{prop:different reps for weighted graded root} and Theorem \ref{thm:invariance}. 
\end{proof}

\begin{notn}
In situations where it becomes helpful to specify the underlying admissible family $F$, the notation $P_{F, [k]}^{\infty}$ will be used rather than $P_{[k]}^{\infty}$.
\end{notn}

\begin{rem}\label{rem:P infinity}
It follows from the definition \eqref{eq:laurent poly enhancement} of the Laurent polynomial weights and from Theorem \ref{thm:convergence} that
\begin{align} \label{P [k] infty}
    P_{F,[k]}^{\infty} = \sum\limits_{x\in \Z^s}F_{\Gamma,k}(x)q^{\eps_{k}(x)}t^{\theta_k(x)},
\end{align}
for any $k\in[k]$. Note for any $F$, setting $t=1$ in $P_{F,[k]}^\infty$ gives a well-defined Laurent $q$-series invariant of $(Y(\Gamma,[k])$.
\end{rem}


\section{The \texorpdfstring{$\widehat{Z}_a(q)$}{} power series}
\label{sec:Zhat_a}

This section starts with a review of the GPPV invariant of negative definite plumbed manifolds, motivating the definition of the admissible family of functions $\widehat F$. In fact, three closely related admissible families are discussed in Section \ref{subsec: Three admissible families}, $\widehat F$, $\widehat F^+$, and $\widehat F^-$. Section \ref{sec: The lattice and Zhat} reformulates the $\widehat Z$ invariant using the lattice cohomology convention for $\spinc$ structures.  
Theorem \ref{thm: Zhathat Zhat} shows that the GPPV invariant is a specialization of the $2$-variable series $P_{\widehat F, [k]}^{\infty}$ at $t=1$, thus establishing Theorem \ref{thm: Zhat intro} stated in the introduction. Additionally, Section \ref{sec: Recovering q series} gives calculations in specific examples.

Let $a\in \delta + 2\Z^s$ be a representative of a $\spinc$ structure $[a]$ on $Y$, using the convention \eqref{eq:GM spinc structures}. Following \cite{GPPV} (see also \cite[Section 4.3]{GM}), consider 
\begin{equation}
\label{eq:Zhat_a}
    \Zhat_a(q) := q^{- \frac{3s + \sum_v m_v}{4}}  \cdot 
    v.p.
    \oint\limits_{\vert z_v \vert =1}  \prod_{v\in \V(\Gamma)} \frac{dz_{v}}{2\pi i z_{v}}\left(z_{v}-\frac{1}{z_{v}}\right)^{2-\delta_v} \cdot \Theta_a^{-M}(z),
\end{equation}
where 
\begin{equation}
\label{eq:theta function}
\Theta_a^{-M}(z) := \sum_{\ell \in a + 2M\Z^s} q^{-\frac{\ell^t M^{-1}\ell}{4}} \prod_{v\in \V(\Gamma)} z_v^{\ell_v}.
\end{equation}
In \eqref{eq:Zhat_a}, $v.p.$ denotes the principal value, that is the average of the integrals over $\vert z_v \vert = 1 +\varepsilon$ and over $\vert z_v \vert = 1 -\varepsilon$, for small $\varepsilon>0$. Note that, since $M$ is negative definite, for each power of $q$ the expression (\ref{eq:theta function}) for
$\Theta_a^{-M}(z)$ is a Laurent polynomial in the variables $\{z_v\}_{v\in \V}$, and the exponent $\left(-\ell^t M^{-1} \ell\right)/4$ as $\ell$ varies is bounded below. It is clear from the definition that $\Zhat_a(q)$ is independent of the choice of representative $a\in [a]$. 

We begin by rewriting \eqref{eq:Zhat_a} as
\begin{align}\label{eq: rewrite}
    \Zhat_a(q)     
    = 
    q^{- \frac{3s+ \sum_{v} m_v}{4}}\sum\limits_{\ell\in a+2M\Z^s} 
    \left[
    \prod\limits_{v}v.p. \frac{1}{2\pi i}\oint\limits_{|z_{v}|=1}\frac{1}{z_{v}}\left(z_{v}-\frac{1}{z_{v}}\right)^{2-\delta_v}z_{v}^{\ell_{v}}dz_{v}
    \right]
    q^{-\frac{\ell^t M^{-1}\ell}{4}}.
\end{align}
Our analysis of the $\widehat Z$ invariant, in particular Proposition \ref{prop:Zhat_k in terms of F} and Theorem \ref{thm: Zhathat Zhat} below, depends on the properties of the coefficient given by the expression in the square brackets in equation \eqref{eq: rewrite}.  Thus we start by rewriting it in more concrete terms. We note that this preliminary analysis, leading to Definition \ref{def:Fhat}, amounts to taking a detailed look at the coefficients denoted $F_{\vec{\ell}}$ in \cite[equation (43)]{GM}.

To compute the integral in equation \eqref{eq: rewrite}, write $\left(z_{v}-z_{v}^{-1}\right)^{2-\delta_v}$ as a Laurent series $E_v^-$ in $z_v$ for the integral over $\lr{z_v} = 1 - \eps$, and as a Laurent series $E_v^+$ in $z_v^{-1}$ for the integral over $\lr{z_v} = 1+ \eps$. Then for $\ell_v\in \Z$, 
\begin{equation*}
\label{eq:principal value with residues}
v.p. \frac{1}{2\pi i}\oint\limits_{|z_{v}|=1}\frac{1}{z_{v}}\left(z_{v}-\frac{1}{z_{v}}\right)^{2-\delta_v}z_{v}^{\ell_{v}}dz_{v}
=
\frac{1}{2}\left[ 
\Res\left(z_v^{\ell_v - 1} E_v^-,0\right) 
+ 
\Res\left(z_v^{-\ell_v - 1}E_v^+(z_v^{-1}),0\right) 
\right].
\end{equation*}
Note that $\Res\left(z_v^{\ell_v - 1} E_v^-,0\right) $ and $\Res\left(z_v^{-\ell_v - 1}E_v^+(z_v^{-1}),0\right) $ equal the coefficient of $z^{-\ell_v}$ in $E_v^-$ and $E_v^+$, respectively.  We will now identify the Laurent series $E_v^-$ and $E_v^+$ more explicitly.

When $\delta_v\leq 2$, the exponent in $(z_v- z_v^{-1})^{2-\delta_v}$ is non-negative and $E_v^+ =E_v^-$ = $(z_v- z_v^{-1})^{2-\delta_v}$ is a Laurent polynomial. In particular, if $\delta_v\leq 2$ for all vertices $v$, then $\Zhat_a(q)$ is a Laurent polynomial with integer coefficients. More generally, coefficients of $\Zhat_a(q)$ are in $2^{-c}\cdot \Z$ where $c$ is the number of vertices of degree at least 3. 

We now describe the Laurent series expansions $E_v^{\pm}$ of $(z_v - z_v^{-1})^{2-\delta_v}$ when $\delta_v\geq 3$. Fix $n\geq 3$. For $\lr{z}<1$, using the expansion $
\left( z - \frac{1}{z}  \right)^{-1}
= 
\frac{-z}{1-z^2}
=
-\sum\limits_{i\geq 0}z^{2i+1} $,
we can write
\begin{equation}
\label{eq:expansion <1}
\left(z - \frac{1}{z}\right)^{2-n} = \left( -\sum\limits_{i\geq 0} z^{2i+1}\right)^{n -2}.
\end{equation}
For $\lr{z}>1$, the expansion $\left( z - \frac{1}{z} \right)^{-1}
= 
\frac{z^{-1}}{1-z^{-2}}
= 
\sum\limits_{i\geq 0}z_v^{-(2i+1)}$ gives 
\begin{equation}
    \label{eq:expansion >1}
    \left(z - \frac{1}{z}\right)^{2-n} = \left( \sum_{i\geq 0} z^{-(2i+1)} \right)^{n-2}.
\end{equation}
Then $E_v^-$ and $E_v^+$ are given by substituting $z=z_v$, $n=\delta_v$ into the right-hand side of \eqref{eq:expansion <1} and \eqref{eq:expansion >1}, respectively. We summarize the discussion so far: the expression in square brackets in equation \eqref{eq: rewrite}
equals the product over $v\in \V(\Gamma)$ of the average of the  coefficients of $z^{-\ell_v}$ in 
\eqref{eq:expansion <1}, \eqref{eq:expansion >1}. 

We now define a family of functions $\Fhat = \{\Fhat_n\}_{n\geq 0}$ which record the coefficients in the average of the two expansions. In Proposition \ref{prop:Fhats are admissible} we show this family is admissible. 

\begin{defn}
\label{def:Fhat}
Consider the following family of functions $\{ \Fhat_n : \Z \to \Q \}_{n\in \Z_+}$.
For $0\leq n \leq 2$, set $\Fhat_n(r)$ to be the coefficient of $z^{-r}$ in $(z-z^{-1})^{2-n}$. For $n\geq 3$, $\Fhat_n(r)$ is defined to be the average of the coefficients of $z^{-r}$ in equations \eqref{eq:expansion <1} and \eqref{eq:expansion >1}.
\end{defn}

Note that $\Fhat_n$ takes values in $\Z$ for $0\leq n\leq 2$ and in $\frac{1}{2}\Z$ for $n \geq 3$. Thus $\Fhat_{\Gamma,k}$, defined in equation \eqref{eq:F_Gamma,k}, takes values in $2^{-c} \Z$ where $c$ is the number of vertices of degree at least $3$. Although an explicit formula for $\Fhat_n$, $n\geq 3$, will not be used in this paper, for the reader's convenience we record it in equation \eqref{eq:Fhat formula}.
\begin{align}
\label{eq:Fhat formula}
    \Fhat_{n}(r) = 
    \begin{cases}
    \frac{1}{2}\sgn(r)^{n}\begin{pmatrix}
    \frac{n+|r|}{2}-2\\
    n-3
    \end{pmatrix}&\text{ if }|r|\geq n-2\text{ and } r\equiv n\bmod 2\\
    0&\text{ otherwise.}
    \end{cases}
\end{align}
Here $\sgn(r) \in \{-1, 1\}$ denotes the sign of $r$.

\subsection{Three admissible families} \label{subsec: Three admissible families}

In this section we introduce families $\Fhat^{+}$, $\Fhat^{-}$, closely related to $\Fhat$, and show that they are all admissible.

\begin{defn}
\label{def:Fhatpm}
For $r\in \Z$ and $0\leq n \leq 2$, set $\Fhat^+_n(r) = \Fhat^-_n(r) = \Fhat_n(r)$ to be the coefficient of $z^{-r}$ in $(z-z^{-1})^{2-n}$. For $n\geq 3$, $\Fhat_n^{-}(r)$ and $\Fhat_n^+(r)$ are defined to be the coefficient of $z^{-r}$ in \eqref{eq:expansion <1} and \eqref{eq:expansion >1}, respectively. 
\end{defn}

The following general observation is used in the proof of the proposition below. If $F^1, F^2, \ldots, F^m$ are admissible families valued in a field of characteristic zero, then the family $\av(F^1, \ldots, F^m)$ given by the average
\begin{equation}
\label{eq:averaging}
\left(\av(F^1, \ldots, F^m) \right)_n = \frac{1}{m} \sum_{i=1}^m F^1_n  + \cdots + F^m_n 
\end{equation}
is again admissible.

\begin{prop}
\label{prop:Fhats are admissible}
The families $\Fhat^+$, $\Fhat^-$, and $\Fhat$ are admissible. 
\end{prop}
\begin{proof}
Property \ref{item:A1} is straightforward to verify. To show property \ref{item:A2} for $\Fhat^+$, note that 
\[
(z-z^{-1}) \sum_{i\geq 0} z^{-(2i+1)} = 1.
\]
Therefore
\[
z \left(  \sum_{i\geq 0} z^{-(2i+1)} \right)^{n-2} - z^{-1} \left( \sum_{i\geq 0} z^{-(2i+1)} \right)^{n-2} =  \left(\sum_{i\geq 0} z^{-(2i+1)} \right)^{n-3},
\]
which demonstrates \ref{item:A2}. Alternatively, \ref{item:A2} may also be seen from a binomial coefficient identity, using an explicit formula for $\widehat F^{\pm}$, analogous to \eqref{eq:Fhat formula}. The calculation for $\Fhat^-$ is similar. Finally, note that $\Fhat$ is the average of $\Fhat^+$ and $\Fhat^-$  and is therefore admissible.
\end{proof}

\subsection{The lattice and \texorpdfstring{$\widehat{Z}$}{}} \label{sec: The lattice and Zhat}

In this section we reformulate $\Zhat$ as a sum of  contributions of the associated function $\Fhat_{\Gamma,k}$ (see equation \eqref{eq:F_Gamma,k} and Remark \ref{Gamma remark}) over lattice points, using the lattice cohomology identification of $\spinc$ structures. 

As a first step, we reparameterize definition \eqref{eq:Zhat_a} in the following way. Every $\ell \in a+2M\Z^s$ can be written in the form $\ell = a + 2M x$ for a unique $x\in \Z^s$. Then 
\[
    \frac{\ell^t M^{-1}  \ell}{4} = \frac{a^2}{4} +   a\cdot x + \langle x,x\rangle
    = \frac{a^2}{4} - 2\chi^{}_{a}(x),
\]
using the notation of Remark \ref{rem: k^2 meaning} and \eqref{eq:chi}. Compare with \cite[Equation (46)]{GM}.  Thus we get 
\begin{equation}
    \label{eq:Zhat_a reformulation}
    \Zhat_a(q) = 
    q^{- \frac{a^2 + 3s + \sum_v m_v}{4}} v.p.
    \oint\limits_{|z_{v}|=1}\prod_{v}\frac{dz_{v}}{2\pi i z_{v}}\left(z_{v}-\frac{1}{z_{v}}\right)^{2-\delta_v}
    \left(\sum\limits_{x\in \mathbb{Z}^{s}}q^{2\chi^{}_{ a}(x)}\prod\limits_{v}z_{v}^{(a+2Mx)_{v}}\right).
\end{equation}

We now move on to the main goal of this section. Recall from Section \ref{sec:identification of spinc structures} that graded roots and lattice cohomology use a different identification of $\spinc$ structures than the $\Zhat$ invariant. The translation between these two identifications is given in equation \eqref{eq:translating between spinc identifications}. Given $k\in m+2\Z^s$, let $a=k-Mu \in \delta + 2\Z^s$ denote the corresponding $\spinc$ representative, and set
\begin{equation*}
    \ZhatL_k(q) := \Zhat_a(q).
\end{equation*}

Recall Notation \ref{notn} for ${\normalization_k}$ and $\eps_k(x)$ in the following statement.

\begin{prop}
\label{prop:Zhat_k in terms of F}
For $k\in m + 2\Z^s$, we have 
\[
\ZhatL_k(q) = \sum_{x\in \Z^s} \Fhat_{\Gamma,k}(x) q^{\eps_k(x)}. 
\]
\end{prop}
\begin{proof}
Note that $2\chi^{}_{a}(x) = 2\chi^{}_{k}(x) + \langle x,u\rangle$ for all $x\in \Z^s$, so equation \eqref{eq:Zhat_a reformulation}  with $a=k-Mu$ can be written as 
\[
\Zhat_a(q) = 
    q^{\normalization_k} v.p.
    \oint\limits_{|z_{v}|=1}\prod_{v}\frac{dz_{v}}{2\pi i z_{v}}\left(z_{v}-\frac{1}{z_{v}}\right)^{2-\delta_v}
    \left(\sum\limits_{x\in \mathbb{Z}^{s}}q^{2\chi^{}_{ k}(x) + \langle x,u\rangle}\prod\limits_{v}z_{v}^{(2Mx+k-Mu)_{v}}\right).
\]
From the above equation and the discussion preceding Section \ref{subsec: Three admissible families}, we see that for every $j\in \Z$,  the coefficient of $q^{\normalization_k + j}$ in $\Zhat_k^{\lattice}$ is equal to
\[
\sum\limits_{\substack{x\in \Z^s \\ 2\chi^{}_{k}(x)+\langle x,u\rangle=j}}\Fhat_{\Gamma, k}(x),
\]
which verifies the desired equality.
\end{proof}

\subsection{Recovering the \texorpdfstring{$q$}{q}-series} \label{sec: Recovering q series}

In this section we show that, when the admissible family is $\Fhat$, the two-variable series specializes to $\Zhat_k^{\lattice}(q)$ by setting $t=1$. Calculations for $S^3$ and $\Sigma(2,7,15)$ are presented further below.

\begin{defn}
Fix a negative definite plumbing $\Gamma$ and a $\spinc$ structure $[k]$. Define
\begin{align}
    \Zhathat_{[k]}(q,t) := P^{\infty}_{\widehat{F},[k]}
\end{align}
which, as we recall from Theorem \ref{thm:convergence}, is an invariant of $(Y(\Gamma), [k])$. 
\end{defn}

\begin{thm} \label{thm: Zhathat Zhat}
With the above notation,
\[
\Zhathat_{[k]}(q,1) = \Zhat_k^{\lattice}(q).
\]
\end{thm}

\begin{proof}
Fix $j\in \Z$. Using the notation in the proof of Theorem \ref{thm:convergence}, the coefficient of $q^{\normalization_k + j}$ in $\Zhathat_{[k]}(q,1)$ is equal to 
\[
\sum_{x\in \partial\til{\sublevel}_j} \Fhat_{\Gamma,k}(x),
\]
which by Proposition \ref{prop:Zhat_k in terms of F} equals the coefficient of $q^{\normalization_k + j}$ in $\Zhat_k^{\lattice}(q)$.
\end{proof}

\begin{exmp}
\label{ex:wgroot -p framed unknot}
Let $\Gamma$ denote the plumbing tree consisting of a single $-p$ framed vertex, with $p\geq 1$. The associated $3$-manifold $Y(\Gamma)$ is the lens space $L(p,1)$, which has $p$ $\spinc$ structures. We illustrate the calculation of the weighted graded root for one $\spinc$ structure; the calculations for the other $\spinc$ structures are analogous. Let $k = -p \in -p + 2\Z$ denote a representative for the $\spinc$ structure $[-p]$. For $x\in \Z$, we have
\[
\chi^{}_k(x) = \frac{p}{2} (x^2+x), \ \ \eps_k(x) =  - \frac{3-p}{4} + px^2, \ \ \theta_k(x) = -px, \ \ \Fhat_{\Gamma,k}(x) = \Fhat_0(-2px).
\]
From the formula \eqref{eq:F_1 and F_0}, we see that $\Fhat_{\Gamma,k}(x) = 0$ unless $x = 0$ or $x= \pm \frac{1}{p}$; this latter case occurs only when $p=1$. The weighted graded root and $\Zhathat_{[k]}(q,t)$ are shown in Figure \ref{fig:wgroot -p}.
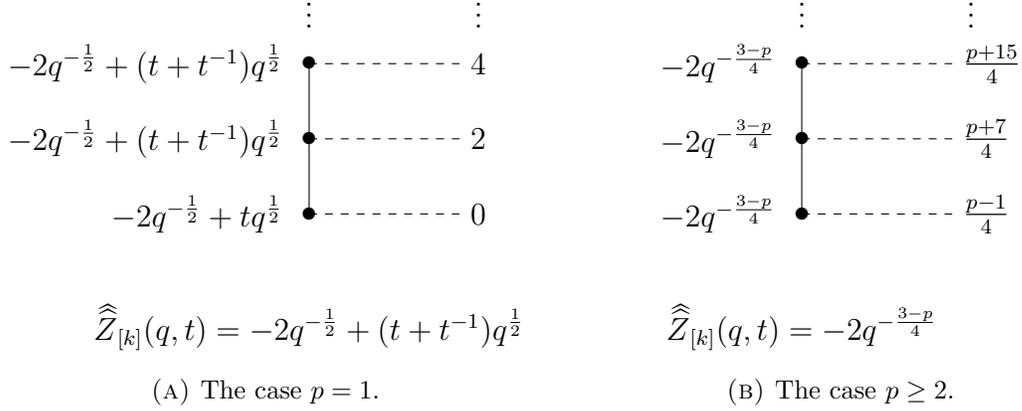
\begin{figure}
\centering
\subcaptionbox{The case $p=1$. \label{fig:wgroot -1}}[.45 \linewidth]{
{   
\begin{tikzpicture}
\node (0) at (0,0) {$\bullet$};
\node (1) at (0,1) {$\bullet$};
\node (2) at (0,2) {$\bullet$};
\node at (0,2.75) {$\vdots$};
\node[right] at (2.05,2.75) {$\vdots$};
\draw (0,0) -- (0,2);

\draw[dashed] (0,0) --++ (2,0);
\draw[dashed] (0,1) --++ (2,0);
\draw[dashed] (0,2) --++ (2,0);

\node[right] at (2,0) {$0$};
\node[right] at (2,1) {$2$};
\node[right] at (2,2) {$4$};

\node[left] at (-.2,0) {$-2q^{-\frac{1}{2}} + tq^{\frac{1}{2}}$};
\node[left] at (-.2,1) {$-2q^{-\frac{1}{2}} + (t+t^{-1})q^{\frac{1}{2}}$};
\node[left] at (-.2,2) {$-2q^{-\frac{1}{2}} + (t+t^{-1})q^{\frac{1}{2}}$};

\node at (0,-1.5) {$\Zhathat_{[k]}(q,t) = -2q^{-\frac{1}{2}} + (t + t^{-1})q^{\frac{1}{2}}$};
\end{tikzpicture}
}
}
\subcaptionbox{The case $p\geq 2$.}[.45 \linewidth]{
\begin{tikzpicture}
\node (0) at (0,0) {$\bullet$};
\node (1) at (0,1) {$\bullet$};
\node (2) at (0,2) {$\bullet$};
\node at (0,2.75) {$\vdots$};
\node[right] at (2.05,2.75) {$\vdots$};
\draw (0,0) -- (0,2);

\draw[dashed] (0,0) --++ (2,0);
\draw[dashed] (0,1) --++ (2,0);
\draw[dashed] (0,2) --++ (2,0);

\node[right] at (2,0) {$\frac{p-1}{4}$};
\node[right] at (2,1) {$\frac{p+7}{4}$};
\node[right] at (2,2) {$\frac{p+15}{4}$};

\node[left] at (-.2,0) {$-2q^{-\frac{3-p}{4}}$};
\node[left] at (-.2,1) {$-2q^{-\frac{3-p}{4}}$};
\node[left] at (-.2,2) {$-2q^{-\frac{3-p}{4}}$};

\node at (0,-1.5) {$\Zhathat_{[k]}(q,t) = -2q^{-\frac{3-p}{4}}$};
\end{tikzpicture}
}
\caption{
The weighted graded root and two-variable series $\Zhathat_{[k]}(q,t)$ corresponding to the admissible family $\Fhat$ for the one-vertex plumbing $\Gamma =\genfrac{}{}{0pt}{}{\mathlarger{\bullet}}{-p}$ with $\spinc$ structure $[-p]$. Recall from grading conventions \ref{grading conv} that the numbers to the right of the graded root denote the gradings of $HF^{+}(-L(p,1),[-p])$.  }
\label{fig:wgroot -p}
\end{figure}

\end{exmp}
In particular, specializing $\Zhathat_{[k]}(q,t)$ at $t=1$ in the case $p=1$ recovers the $\Zhat$ invariant for $S^3$. The calculation above shows that the $2$-variable series $\Zhathat$ introduced in this paper is different from the conjectured Poincar\'{e} series of the BPS homology \cite[Equation (6.80)]{GPV}, \cite[Equation (18)]{GM}.

\begin{exmp} \label{2 7 15 ex involution0} 
Consider the Brieskorn sphere $\Sigma(2,7,15)$, which can be represented as a negative definite plumbing, as shown in Figure \ref{fig:plumbing and link ex}. Since $\Sigma(2,7,15)$ is an integer homology sphere, we denote its unique $\spinc$ structure by $\mathfrak{s}_0$.

One can compute, cf. \cite[Section 4.6]{GM}:
\begin{align*}
    \Zhat_{\mathfrak{s}_{0}}(q) &=
    q^{1739/840}\sum\limits_{n=0}^{\infty}\left[q^{(61+420n)^2/840}+q^{(149+420n)^2/840}+q^{(299+420n)^2/840}+q^{(331+420n)^2/840}\right]\\
    &-q^{1739/840}\sum_{n=0}^{\infty}\left[q^{(89+420n)^2/840}+q^{(121+420n)^2/840}+q^{(271+420n)^2/840}+q^{(359+420n)^2/840}\right]\\
    &=q^{13/2} - q^{23/2} - q^{39/2}+q^{57/2} - q^{179/2} + q^{217/2} + q^{265/2} - q^{311/2} +\cdots 
\end{align*}
The beginning of the weighted graded root (a result of a computer calculation \cite{program}) for $\Sigma(2,7,15)$ is shown in Figure \ref{Weighted_graded_root}; additional weights are given in the table below. 

\begin{table}[H]
\begin{tabular}{|c|c|}
\hline
Weight                                                                                                                                           & Grading \\[.2cm]\hline
\rule{0pt}{3ex}$\frac{1}{2}(t + t^{-1})q^{\frac{13}{2}} - q^{\frac{23}{2}} - q^{\frac{39}{2}}$                                 & 20      \\[.2cm] \hline
\rule{0pt}{3ex} $\frac{1}{2}(t + t^{-1})q^{\frac{13}{2}} - q^{\frac{23}{2}} - q^{\frac{39}{2}} + \frac{1}{2}(t+t^{-1})q^{\frac{57}{2}}$ & 28      \\[.2cm] \hline
\end{tabular}
\end{table}
\noindent In particular, setting $t=1$ in the weight at grading $28$, one can see the first few terms of $\Zhat_{\mathfrak{s}_{0}}(q)$ as result of stabilization, which is a consequence of Theorems \ref{thm:convergence}, \ref{thm: Zhathat Zhat}.
\end{exmp}

\section{\texorpdfstring{$\Spinc$}{Spinc} conjugation}\label{section:conj}
In this section we study the behavior of $\Zhathat$ and weighted graded roots under $\spinc$ conjugation. Under the identification (\ref{eq:LC spinc structures}), conjugation is given by the map $[k]\to [-k]$. Both $\Zhat$ and lattice cohomology, in particular graded roots, are invariant under conjugation. However, when considering our new theory of weighted graded roots, a different, more refined, story emerges which we now describe. 

Let $F$ be an $\Ring$-valued admissible family. Consider the following property. 
\begin{equation}
\label{eq:A3}
    F_n(-r) = (-1)^n F_n(r) \text{ for all } n\geq 0 \text{ and } r\in \Z. \tag{A3} 
\end{equation}

\begin{prop}\label{t inversion}
If $F$ is an admissible family which satisfies property \eqref{eq:A3}, then $P_{F,[k]}^{\infty}(q,t) = P_{F,[-k]}^{\infty}(q,t^{-1})$ for all $\spinc$ structures $[k]$. 
\end{prop}
\begin{proof}
Note that $k':= -k + 2Mu$ is a representative for $[-k]$. We will show that 
\begin{equation} \label{eq: Zhathat conj proof}
F_{\Gamma,k}(x) q^{\eps_k(x)}t^{\theta_k(x)} = F_{\Gamma,k'}(-x)q^{\eps_{k'}(-x)}(t^{-1})^{\theta_{k'}(-x)}
\end{equation}
for all $x\in \Z^s$, and the claim follows. First,
\[
\Delta_{k'} = \Delta_{k},\hspace{2em} 2\chi_k(x)  = 2\chi^{}_{k'}(-x) - 2 \langle x,u\rangle,
\]
so $\eps_{k}(x) = \eps_{k'}(-x)$ for all $x\in \Z^s$. Next, $2Mx + k - Mu = -\left(2M(-x) +k' - Mu \right)$, so
\begin{align*}
    F_{\Gamma,k}(x) = (-1)^{\sum_v \delta_v} F_{\Gamma,k'}(-x) = F_{\Gamma,k'}(-x),
\end{align*}
where the first equality follows from property \eqref{eq:A3} and the second is due to the sum of degrees in any graph being even. Lastly, 
\begin{align*}
    \tnormalization_{k'}=\frac{(-k+2Mu)\cdot u-\langle u, u\rangle}{2} = -\frac{k\cdot u-\langle u, u\rangle }{2} = -\tnormalization_k
\end{align*}
So $\theta_{k'}(-x) = -\tnormalization_k+\langle -x,u\rangle = -\theta_k(x)$.
\end{proof}

Now note that $\Fhat$, introduced in Definition \ref{def:Fhat}, satisfies property \eqref{eq:A3}.\footnote{Although it will not be used, we note that $\Fhat^{\pm}$ from Definition \ref{def:Fhatpm} do not satisfy \eqref{eq:A3}.}
\begin{cor} \label{prop: Zhathat conjugation}
$\Zhathat_{[k]}(q, t) = \Zhathat_{[-k]}(q, t^{-1})$. In particular, setting $t = 1$ recovers the conjugation invariance of $\Zhat$. 
\end{cor}

We now turn to weighted graded roots and illustrate, via two examples, some interesting behavior under $\spinc$ conjugation. First, we briefly recall how graded roots transform under conjugation and describe the corresponding story in Heegaard Floer homology. 

Given a negative definite plumbing $\Gamma$ and a $\spinc$ structure $[k]$, the map $\Z^s\to \Z^s$, sending $x$ to $-x$
induces an isomorphism $(R_{k}, \chi^{}_{k})\cong (R_{-k}, \chi^{}_{-k})$ since $\chi^{}_{k}(x) = \chi^{}_{-k}(-x)$. Similarly, in Heegaard Floer homology, for any closed oriented 3-manifold $Y$ and $\spinc$ structure $\mathfrak{s}$, there is an isomorphism $HF^+(Y,\mathfrak{s})\cong HF^+(Y, \bar{\mathfrak{s}})$, where $\bar{\mathfrak{s}}$ is the conjugate of $\mathfrak{s}$; see \cite[Theorem 2.4]{OS-applications}. 

Moreover, for a self-conjugate $\spinc$ structure we get an involution on the graded root and on Heegaard Floer homology. The involution on the graded root is induced by the map
\begin{align*}
    \Z^s\to \Z^s,\hspace{1em} x\mapsto -x-M^{-1}k.
\end{align*}
Note here $M^{-1}k\in \Z^s$ since $[k]=[-k]$. For Heegaard Floer homology, the involution
\begin{align*}
    \iota: HF^+(Y, \mathfrak{s}) \to HF^+(Y, \mathfrak{s})
\end{align*}
comes from a chain map obtained by considering what happens when a pointed Heegaard diagram $(\Sigma, \alpha, \beta, z)$ representing $Y$ is replaced with $(-\Sigma, \beta, \alpha, z)$. The involution $\iota$ is at the foundation of involutive Heegaard Floer homology, an extension of Heegaard Floer homology due to Hendricks-Manolescu  \cite{involutive}. Note, involutive Heegaard Floer homology is currently only defined over $\F = \Z/2\Z$. So when discussing the involution $\iota$, we will assume we are working with $\F$ coefficients.

For $\Gamma$ an almost rational plumbing, Dai-Manolescu show that the two involutions described above are identified under the isomorphism given in Theorem \ref{thm: LC-HF iso} (see \cite[Theorem 3.1]{DM}). Furthermore, they show that the graded root is symmetric about the infinite stem and the involution is the reflection about the infinite stem.  

\begin{exmp} \label{2 7 15 ex involution} 
Consider again the Brieskorn sphere $\Sigma(2,7,15)$. Note, the plumbing given in Figure \ref{fig:plumbing and link ex} describing $\Sigma(2,7, 15)$ is almost rational. Also, since $\Sigma(2, 7, 15)$ only has one $\spinc$ structure, $\mathfrak{s}_{0}$, it is self-conjugate by default. Hence, the corresponding graded root is symmetric about the infinite stem and the involution is the reflection. However, as seen in Figure \ref{Weighted_graded_root}, the weighted graded root is no longer symmetric and the involution does not preserve all of the weights. There is a node at grading level $6$ which has weight $\frac{1}{2}tq^{13/2}$, whereas the node on the opposite side of the infinite stem has weight $0$. The reason for this symmetry breaking is a result of the failure of $\Fhat_{\Gamma, k}(-x-M^{-1}k)q^{\eps_{k}(-x-M^{-1}k)}t^{\theta_k(-x-M^{-1}k)}$ to equal $\Fhat_{\Gamma, k}(x)q^{\eps_{k}(x)}t^{\theta_k(x)}$.
\end{exmp}

The following example shows that, unlike $\Zhat$ and graded roots, the weighted graded root can distinguish conjugate $\spinc$ structures. Moreover, it exhibits a new phenomenon different from that in Corollary \ref{prop: Zhathat conjugation}.

\begin{exmp} \label{ex: asymmetric conjugation}
Let $\Gamma$ be the plumbing pictured below:
\begin{center}
\begin{tikzpicture}[scale=.7]
\node (m) at (0,0) {$\bullet$};
\node (r) at (2,0) {$\bullet$};
\node (l) at (-2,0) {$\bullet$};
\node (b) at (0,-2) {$\bullet$};
\node (t) at (0,2) {$\bullet$};

\draw (-2,0) --++ (4,0);
\draw (0,2) --++ (0,-4);

\node[above] at (-.5,0) {$-1$};
\node[above] at (r) {$-7$};
\node[above] at (l) {$-11$};
\node[left] at (b) {$-10$};
\node[left] at (t) {$-3$};

\end{tikzpicture}
\end{center}
Order the vertices so that $v_{1}, v_{2}, v_{3}, v_{4}, v_{5}$ correspond to the vertices with weights $-1$, $-7$,$-10$,\\$-11$,$-3$. Let $k=(-5,5,8,9, 1)$. Consider the $\spinc$ structure $[k]$ and its conjugate $[-k]$.

\begin{figure}
\centering
{\begin{tikzpicture}
    \node at (-3,0) {$\bullet$};
    \node at (-3,1) {$\bullet$};
    \node at (-3,2) {$\bullet$};
    \node at (-3,3) {$\bullet$};
    \node at (-3,4) {$\bullet$};
    \node at (-3,5) {$\bullet$};
    \node at (-3,6) {$\bullet$};
    \node at (-3,7) {$\bullet$};
    \node at (-3,8) {$\bullet$};
    \node at (-3,9) {$\bullet$};
    
    \node at (3,0) {$\bullet$};
    \node at (3,1) {$\bullet$};
    \node at (3,2) {$\bullet$};
    \node at (3,3) {$\bullet$};
    \node at (3,4) {$\bullet$};
    \node at (3,5) {$\bullet$};
    \node at (3,6) {$\bullet$};
    \node at (3,7) {$\bullet$};
    \node at (3,8) {$\bullet$};
    \node at (3,9) {$\bullet$};
    
    \draw (-3,0)--(-3, 9);
    \draw (3,0)--(3, 9);
    
    \node[right] at (7,0) {$\frac{-570}{769}$};
    \node[right] at (7,1) {$\frac{968}{769}$};
    \node[right] at (7,2) {$\frac{2506}{769}$};
    \node[right] at (7,3) {$\frac{4044}{769}$};
    \node[right] at (7,4) {$\frac{5582}{769}$};
    \node[right] at (7,5) {$\frac{7120}{769}$};
    \node[right] at (7,6) {$\frac{8658}{769}$};
    \node[right] at (7,7) {$\frac{10196}{769}$};
    \node[right] at (7,8) {$\frac{11734}{769}$};
    \node[right] at (7,9) {$\frac{13272}{769}$};
    
    \draw[dashed] (-3,0)--(7,0);
    \draw[dashed] (-3,1)--(7,1);
    \draw[dashed] (-3,2)--(7,2);
    \draw[dashed] (-3,3)--(7,3);
    \draw[dashed] (-3,4)--(7,4);
    \draw[dashed] (-3,5)--(7,5);
    \draw[dashed] (-3,6)--(7,6);
    \draw[dashed] (-3,7)--(7,7);
    \draw[dashed] (-3,8)--(7,8);
    \draw[dashed] (-3,9)--(7,9);
    
    \node[left] at (3,.2) {$0$};
    \node[left] at (3,1.2) {$0$};
    \node[left] at (3,2.2) {$0$};
    \node[left] at (3,3.2) {$0$};
    \node[left] at (3,4.3) {$\frac{1}{2}t^3q^{\frac{15009}{1538}}$};
    \node[left] at (3,5.3) {$\frac{1}{2}t^3q^{\frac{15009}{1538}}$};
    \node[left] at (3,6.3) {$\frac{1}{2}t^3q^{\frac{15009}{1538}}$};
    \node[left] at (3,7.3) {$q^{\frac{19623}{1538}} + \frac{1}{2}t^3q^{\frac{15009}{1538}}$};
    \node[left] at (3,8.3) {$q^{\frac{19623}{1538}} + \frac{1}{2}t^3q^{\frac{15009}{1538}}$};
    \node[left] at (3,9.3) {$q^{\frac{19623}{1538}} + \frac{1}{2}t^3q^{\frac{15009}{1538}}$};
    
    \node[left] at (-3,.2) {$0$};
    \node[left] at (-3,1.2) {$0$};
    \node[left] at (-3,2.2) {$0$};
    \node[left] at (-3,3.2) {$0$};
    \node[left] at (-3,4.2) {$0$};
    \node[left] at (-3,5.2) {$0$};
    \node[left] at (-3,6.2) {$0$};
    \node[left] at (-3,7.3) {$q^{\frac{19623}{1538}}
    +\frac{1}{2}t^{-3}q^{\frac{15009}{1538}}$};
    \node[left] at (-3,8.3) {$q^{\frac{19623}{1538}}
    +\frac{1}{2}t^{-3}q^{\frac{15009}{1538}}$};
    \node[left] at (-3,9.3) {$q^{\frac{19623}{1538}}
    +\frac{1}{2}t^{-3}q^{\frac{15009}{1538}}$};
    
    \node at (-3,10) {$\vdots$};
    \node at (3,10) {$\vdots$};
    \node at (7.5,10) {$\vdots$};
    
    \node at (-3,-1) {$(R_k, \chi^{}_{k}, P_{\Fhat,k})$};
    \node at (3,-1) {$(R_{-k}, \chi^{}_{-k}, P_{\Fhat,-k})$};
    \end{tikzpicture}}
    \caption{Note, $d(-Y(\Gamma), [(-5,5,8,9,1)]) = d(-Y(\Gamma), [(5,-5,-8,-9,-1)]) = -\frac{570}{769}$  \label{fig: cross graded root}. }
\end{figure}
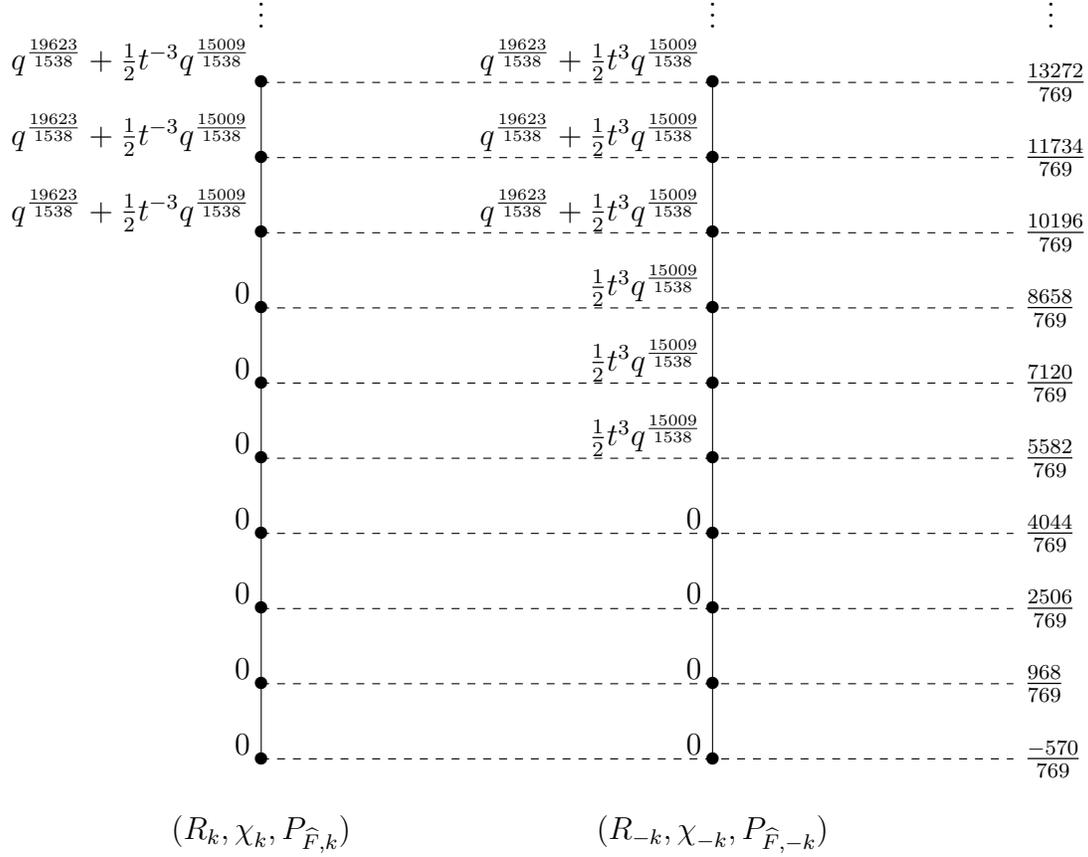
Initial segments of the weighted graded roots (a result of a computer calculation \cite{program})  corresponding to $[k]$ and $[-k]$ are pictured in Figure \ref{fig: cross graded root}.  As discussed in the beginning of this section, the $\Zhat$ invariant and graded roots are invariant under $\spinc$ conjugation. In this example the weighted graded roots not only distinguish $[k]$ and $[-k]$, they do this by more than just inversion of $t$ in all the weights, compare with Corollary \ref{prop: Zhathat conjugation}. 
For example, the node at grading level $\frac{5582}{769}$ for $(R_{[k]}, \chi^{}_{[k]}, P_{\Fhat,[k]})$ is $0$, while the corresponding node for $(R_{[-k]}, \chi^{}_{[-k]}, P_{\Fhat,[-k]})$ is $\frac{1}{2}t^3q^{\frac{15009}{1538}}$.

Note that $\Zhathat(q,t)$ is the limit of the Laurent polynomial weights, whose coefficients stabilize in every bidegree according to Theorem \ref{thm:convergence}. The weighted graded roots carry the unstable information as well; this explains the discrepancy between this example and Corollary \ref{prop: Zhathat conjugation}.
On a more detailed level, the reason for the discrepancy by more than just inversion of $t$ is due to the failure of $\Fhat_{\Gamma, -k}(-x)q^{\eps_{-k}(-x)}(t^{-1})^{\theta_{-k}(-x)}$ to equal $\Fhat_{\Gamma, k}(x)q^{\eps_{k}(x)}t^{\theta_k(x)}$. Equation
\eqref{eq: Zhathat conj proof} in the proof of Proposition \ref{t inversion}, where $k'= -k + 2Mu$,
was sufficient for showing $\Zhathat_{[k]}(q, t)= \Zhathat_{[-k]}(q, t^{-1})$ because the sum is taken over all lattice points $x\in \Z^s$. However, the weights on the nodes of the graded root are sums over lattice points in some connected component of a sublevel set of $\chi^{}_{k}$ for $(R_k, \chi^{}_{k}, P_{\Fhat,k})$, and of $\chi^{}_{-k}$ for $(R_{-k}, \chi^{}_{-k}, P_{\Fhat,-k})$. But the map $x\mapsto -x$ takes the connected components of $\chi^{}_{k}$ sublevel sets to connected components of $\chi^{}_{-k}$ sublevel sets, not connected components of $\chi^{}_{-k + 2Mu}$ sublevel sets.

\end{exmp}


\bibliographystyle{amsalpha}
\bibliography{main}
\nocite{*}
\end{document}